\documentclass[12pt,reqno]{amsart}
\usepackage{amscd,amsmath,amsxtra,amssymb,amsthm,latexsym,amsfonts}
\usepackage[all]{xy}
\usepackage{epic,eepic}
\def\aut{\text{Aut}}
\def\hom{\text{Hom}}

\def\ext{\text{Ext}}

\def\im{\text{Im}}
\def\coker{\text{Coker}}

\def\coh{\text{Coh}}

\def\hilb{\text{Hilb}}

\def\reg{\text{reg}}

\def\3m{M_{3m+1}}
\def\gl{\text{GL}}

\def\GL{\text{GL}}
\def\Coker{\text{Coker}}
\def\Bl{\text{Bl}}

\swapnumbers

\newtheorem{lemma}{Lemma}[section]
\newtheorem{theorem}[lemma]{Theorem}

\newtheorem{corollary}[lemma]{Corollary}
\newtheorem{proposition}[lemma]{Proposition}

\newtheorem{sub}[lemma]{}

\newcommand{\A}{{\mathbb A}}
\newcommand{\Z}{{\mathbb Z}}

\newcommand{\C}{{\mathbb C}}
\newcommand{\Q}{{\mathbb Q}}

\renewcommand{\P}{{\mathbb P}}

\newcommand{\F}{{\mathbb F}}

\newcommand{\ka}{{\mathcal A}}
\newcommand{\kb}{{\mathcal B}}

\newcommand{\ke}{{\mathcal E}}
\newcommand{\kf}{{\mathcal F}}
\newcommand{\kg}{{\mathcal G}}
\newcommand{\kh}{{\mathcal H}}
\newcommand{\ki}{{\mathcal I}}

\newcommand{\kk}{{\mathcal K}}
\newcommand{\kl}{{\mathcal L}}
\newcommand{\km}{{\mathcal M}}
\newcommand{\kn}{{\mathcal N}}
\newcommand{\ko}{{\mathcal O}}

\newcommand{\kq}{{\mathcal Q}}

\newcommand{\kz}{{\mathcal Z}}

\newcommand{\bkf}{{\boldsymbol{\mathcal F}}}
\newcommand{\bkg}{{\boldsymbol{\mathcal G}}}

\newcommand{\bku}{{\boldsymbol{\mathcal U}}}
\newcommand{\bkv}{{\boldsymbol{\mathcal V}}}

\newcommand{\tko}{{\widetilde{\ko}}}

\newcommand{\fraca}{{\mathfrak a}}

\newcommand{\sw}{{\scriptstyle \wedge}}

\CompileMatrices
\SelectTips{eu}{12}
\xyrequire{poly}

\begin{document}
\title[{Sheaves on cubic Space Curves}]{On the Moduli Scheme of stable Sheaves
  supported on cubic Space Curves}

\author[{H.~G. Freiermuth}]{Hans Georg Freiermuth}
\address{Department of Mathematics\\
         Columbia University\\
         2990 Broadway\\
         New York, NY 10027}
\email{freiermuth@math.columbia.edu}


\author[{G. Trautmann}]{G\"unther Trautmann}
\address{Fachbereich Mathematik\\
         Universit\"at Kaiserslautern\\
         Erwin-Schr\"odinger Strasse\\
         D-67663 Kaiserslautern}
\email{trm@mathematik.uni-kl.de}
\date{April 11, 2002}


\begin{abstract}
We investigate the geometry of the Simpson moduli space
$M_{P}(\P_3)$ of stable sheaves with Hilbert polynomial $P(m)=3m+1$.
It consists of two smooth, rational components $M_0$ and $M_1$ of
dimensions $12$ and $13$ intersecting each other transversally
along an $11$-dimensional, smooth, rational subvariety. The
component $M_0$ is isomorphic to the closure of the space of twisted
cubics in the Hilbert scheme $\hilb_P(\P_3)$ and $M_1$ is isomorphic
to the incidence variety of the relative Hilbert scheme of cubic
curves contained in planes. In order to obtain the result and to
classify the sheaves, we characterize $M_P(\P_3)$ as geometric
quotient of a certain matrix parameter space by a non-reductive group.
We also compute the Betti numbers of the Chow groups of the moduli
space.
\end{abstract}

\maketitle
\tableofcontents
\newpage
\section{Introduction}

Generalizing the Maruyama schemes of semi-stable torsionfree coherent
sheaves, C.~Simpson introduced in 1994 coarse projective moduli spaces
$M_P(X)$ for {\sl arbitrary} semi-stable sheaves with fixed Hilbert
polynomial $P$ on a smooth projective variety $X$. If the degree $d$
of the Hilbert polynomial $P$ is strictly less than the dimension of
$X$, all elements of the space $M_P(X)$ are sheaves supported on
proper $d$-dimensional subvarieties of $X$.

One among many reasons for studying these Simpson moduli spaces is
that they often contain structure sheaves of subvarieties of $X$ with
Hilbert polynomial $P$. Thus, certain components of $M_P(X)$ can be
viewed as compactifications of an open part of the corresponding
Hilbert scheme $\hilb_P(X)$.

The structure of the spaces $M_P(\P_2)$ of semi-stable sheaves with
linear polynomial $P(m)=\mu\,m+\chi$ and small multiplicities
$\mu\leq 4$ supported on {\sl plane curves} has been described by
J.~Le Potier \cite{LePot2}.

In this paper, we investigate the geometry of the space
$M_{3m+1}(\P_3)$, one of the first nontrivial examples of a reducible moduli
space of stable sheaves on $\P_3$ supported on {\sl space curves}.

One might view $M=M_{3m+1}(\P_3)$ in some sense as an analogue of the
Hilbert scheme $\hilb_{3m+1}(\P_3)$ containing the twisted cubics,
which has been analyzed by R.~Piene and M.~Schlessinger \cite{P/S}.
The two spaces share many similar properties. They consist for example
both of two smooth rational components which intersect transversally.
Indeed, even more is true: One component $M_0$ of the Simpson moduli
space $M$ is isomorphic to the closure $H_0$ of the space of
twisted cubic curves in $\hilb_{3m+1}(\P_3)$. On the other hand, the
second component $M_1\subset M$ consisting of Cohen-Macaulay modules
on planar cubics is simpler than the ``ghost'' component $H_1$ of the
Hilbert scheme.

The structure sheaves of cubic space curves $\ko_C\in M_0$ specialize
in direction towards the intersection $M_0\cap M_1$ in a natural way
to rank-$1$ Cohen-Macaulay sheaves with support on {\sl singular}
plane cubics.

We give a complete classification of the sheaves in the moduli space
$M$ and of their resolutions. Furthermore, we characterize $M$ as a geometric
quotient of a parameter space $X$ by a non-reductive algebraic group
$G$. The points in $X$ are the matrices that occur in a common
resolution for all the sheaves in $M$. This representation as a
quotient is used to compute an explicit deformation of the sheaves in
$M_0\cap M_1$ into structure sheaves of space cubics, to prove the
transversality of the intersection $M_0\cap M_1$, and to identify $M_0$ as a blow
up. We summarize our result in the following theorem.

\vskip5mm

\begin{theorem}\label{MainThm}
Let $M:=M_{3m+1}(\P_3)$ be the Simpson moduli space of semi-stable
sheaves on $\P_3$ with Hilbert polynomial $3m+1$. Then
\begin{itemize}\itemsep3pt
\item[(1)] $M$ is a fine moduli space representing the corresponding functor.
\item[(2)] All sheaves in $M$ are stable and
  admit a free resolution
$$\xymatrix{0\ar[r] & 2\ko_{\P_3}(-3)\ar[r]^-{B} & \ko_{\P_3}(-1)\oplus
  3\ko_{\P_3}(-2)\ar[r]^-{A}& \ko_{\P_3}\oplus\ko_{\P_3}(-1)\ar[r] &
  \kf \ar[r] & 0}.$$
\item[(3)] $M$ is a projective variety and consists
  of two nonsingular, irreducible, rational
  components $M_0$ and $M_1$ of dimension $12$ and $13$ respectively.
\item[(4)] $M_0$ is isomorphic to the component $H_0$ of the Hilbert
  scheme $\text{\rm Hilb}_{\,3m+1}(\P_3)$ which contains the twisted
  cubic curves. It is also isomorphic to the blow up of a space $N$ of nets
  of quadrics in $\P_3$ along the subvariety $N_1$ of degenerate nets,
  see \cite{E/S} and \ref{nets}.
\item[(5)] The component $M_1$ consists of isomorphism classes of stable
  sheaves $\kf$, supported on plane cubics $C$, which admit a non-split
  extension
  $$\xymatrix{0 \ar[r]&\ko_C\ar[r]&\kf\ar[r]&k_p\ar[r]&0},$$
  where $p\in C$. $M_1$ is isomorphic to the incidence variety of the
  relative Hilbert scheme of cubic curves contained in planes. As a
  subscheme of $M$, the component $M_1$ is characterized by the
  vanishing of the degree $0$ entry of the matrices $A$ in $(2)$.
\item[(6)] The two components $M_0$ and $M_1$ intersect transversally and the
  intersection is smooth, rational of dimension $11$. It is
  isomorphic to the exceptional divisor $\P(\kn_{N_1/N})$ of the blow
  up $M_0\cong\text{Bl}_{N_1}(N)$. The elements of $M_0\cap M_1$ are
  Cohen-Macaulay sheaves
  supported on plane singular cubics with a section vanishing exactly
  in one of the singularities of the curve.
\item[(7)] There is a quasi--affine scheme $X$ acted on by a
  non-reductive algebraic group
  $G$ and a $G$--equivariant morphism $X\xrightarrow{\pi} M$ such that
  \begin{enumerate}
  \item [(i)] the pullback to $X\times \P_3$ of the universal sheaf
    $\bku$ on $M\times\P_3$ has a universal resolution which restricts to the
    resolution of $\bku_{\pi(x)}$ of $(2)$ at any $x$.
  \item [(ii)]$X\xrightarrow{\pi} M$ is a geometric quotient.
  \end{enumerate}
\end{itemize}
\end{theorem}
The description of $M_0$ as a blow up of $N$ is used in section 9 to compute 
the Chow groups of $M$ based on a result of G.~Ellingsrud and
S.A.~Str{\o}mme \cite{E/S} on the Chow ring of the space of nets of
quadrics, see theorem \ref{chowgr}.

\vskip5mm

{\bf Acknowledgements.} The first author wants to thank Michael
Thaddeus for useful discussions and his encouragement.
The second author was partially supported by the
DFG--Schwerpunkt\-programm ``Globale Methoden in der komplexen Geometrie''.

\vskip5mm

{\bf Notations.} We work over a fixed, algebraically closed field $k$ of characteristic
$0$ and denote projective $n$-space by $\P=\P_n=\P V$, where $V$ is an
$(n+1)$-dimensional vector space over $k$. Let $S$ be a Noetherian
(base-)scheme of finite type over $k$. For any (not necessarily reduced
or irreducible) projective variety $X$ with very ample line bundle
$\ko_X(1)$, we denote the two projections
from $S\times_k X$ to $S$ and $X$ by $p$ and $q$ respectively.

Let $\bkf\in\coh(S\times_k X)$, $\kg\in\coh(S)$ and $\kh\in\coh(X)$ be
coherent sheaves. We will use the following abbreviations:
$\kg\boxtimes\kh = p^{\ast}\kg\otimes
q^{\ast}\kh$, $\widetilde{\kg}(m) = \kg\boxtimes \ko_X(m)= p^\ast\kg\otimes q^{\ast}\ko_X(m)$,
$\,\bkf(m) = \bkf \otimes p^{\ast}\ko_X(m)$ and
$\bkf_s=\bkf|_{X_s}$ for the fibers. If $f\,:\,S'\to S$ is a morphism we write $f_X$ for $f\times id  \,:\,
S'\times_k X \to S\times_k X$. Usually, we will use bold letters for
$S$-flat families $\bkf$ on $S\times_k X$. Finally, we denote the support of $\kh\in\coh(X)$ defined by the
annihilator ideal sheaf $\text{Ann}(\kh)$ by $Z(\kh)$.

\vskip10mm

\section{Schemes of stable sheaves}

We recall the theorem of C. Simpson \cite{Sim} on the existence of a
coarse moduli scheme for semi-stable sheaves on
a smooth projective variety $X$. A coherent $\ko_X$--module is
called {\bf purely $d$--dimensional} if its support is purely $d$--dimensional 
and if it has no torsion in dimension $<d$. The Hilbert polynomial of
a purely $d$--dimensional sheaf $\kf$ can be written as
\[
P_{\kf}(m):=\sum\limits_{\nu=0}^d (-1)^\nu\dim_k H^\nu(X,
\kf\otimes\ko_X(m))=\sum_{\nu=0}^d a_\nu(\kf)\binom{m+\nu-1}{\nu}=
\dfrac{a_d(\kf)}{d!}m^d + \ldots
\]
with integers $a_{\nu}(\kf)$. The number $\mu(\kf)=a_d(\kf)$ is positive 
and is called the {\bf multiplicity} of $\kf$. The polynomial
\[
R_\kf (m)=P_\kf(m)/a_d(\kf)
\]
is called the {\bf reduced Hilbert polynomial} of $\kf$. We write $R_\ke\leq R_\kf$
resp. $R_\ke< R_\kf$ if $R_\ke(m)\leq R_\kf(m)$ resp. $R_\ke(m)<R_\kf(m)$ for
large $m$. A coherent sheaf $\kf$ of pure dimension $d$ is called {\bf
  semi-stable} resp. {\bf stable} if for any proper subsheaf $0\neq \ke
\underset{\neq}{\subset} \kf$,
\[
R_\ke\leq R_\kf\quad\text{ resp. }\quad R_\ke<R_\kf.
\]

It is easy to prove that $\kf$ is (semi)-stable if and only if for any proper
quotient sheaf $\kg$ of $\kf$ which is also purely $d$--dimensional,
\[
R_\kf\leq R_\kg\quad\text{ resp. }\quad R_\kf< R_\kg,
\]
see \cite{H/L}, proposition 1.2.6. 

Given a numerical polynomial $P\in \Q[T]$, a contravariant
functor
\[
\km_P(X): \text{(schemes)}\to \text{ (sets) }
\]
is defined as follows. For a scheme $S$, the set $\km_P(X)(S)$ is
the set of all classes $[\bkf]$ of coherent sheaves $\bkf$ on $S\times X$
which are $S$--flat and for which all fibers $\kf_s$, $s\in S$, are semi-stable
with Hilbert polynomial $P_{\kf_s}(m)=P(m)$. Here the class $[\bkf]$ is defined
to be the class of $\bkf$ under the equivalence relation $\bkf\sim\bkf\otimes
p^\ast\kl$ for a line bundle $\kl$ on $S$. For a morphism $T\xrightarrow{f} S$
the map
\[
\km_P(X)(S)\to\km_P(X)(T)
\]
is defined by $[\bkf]\mapsto [f_X^\ast\bkf]$.
\vskip5mm

 \begin{theorem}[Simpson] Let $X$ be a smooth projective variety and $P$ a
      numerical polynomial. Then there exists a projective variety
      $M_P(X)$ which is a coarse moduli space corepresenting the functor
      $\km_P(X)$.\\
      The closed points of $M_P(X)$ parametrize
      the $S$-equivalence classes of semi-stable sheaves on $X$ with fixed Hilbert
      polynomial $P$.\\
      Furthermore, there is an open subset $M_P^s(X)\subset M_P(X)$
      parametrizing the isomorphism classes of stable sheaves.
 \end{theorem}

For proofs, see \cite{Sim},\cite{LePot1},\cite{H/L}. If the integers $a_\nu(\kf)$ in
the Hilbert polynomial are pairwise coprime, then any semi-stable sheaf is already
stable and thus $M_P(X)=M_P^s(X)$. \rm Moreover, in this case there exists a
universal sheaf $\bku$ on $M_P(X)\times X$ such that $M_P(X)$ becomes a fine
moduli space and represents the functor $\km_P(X)$, see \cite{H/L}, corollary
4.6.6. There is also a relative version of Simpson's theorem, see \cite{H/L},
theorem 4.3.7.
\vskip5mm

\begin{proposition}\label{tangen} The tangent space at a stable sheaf $[\kf]$ in
$M_P^s(X)$ is isomorphic to $\ext^{\,1}_{\,\ko_X}(X, \kf, \kf)$ and the germ of
  $M_P^s(X)$ at $[\kf]$ is a universal deformation of $\kf$. If
  $\ext_{\,\ko_X}^{\,2}(X, \kf,\kf)=0$, then $M_P^s(X)$ is smooth at $[\kf]$.
\end{proposition}

A proof can be found in \cite{H/L}, corollary 4.5.2. For the
polynomial $P(m)=3m+1$ we obtain in particular:
\vskip5mm

\begin{proposition}\label{genmod}
a) Any semi-stable sheaf on $\P_3$ with Hilbert polynomial $3m+1$ is
stable and locally Cohen-Macaulay on its support, and 
$M=M_{3m+1}(\P_3)=M^s_{3m+1}(\P_3).$

b) $M$ represents the functor $\km_{3m+1}(\P_3)$ and there is a
  universal sheaf\, $\bku$ over $M\times \P_3$.
\end{proposition}
\vskip10mm

\section{Stable sheaves supported on cubics}

\rm The Hilbert Scheme $H:=\hilb_{3m+1}(\P_3)$
consists of two smooth, rational
components. One of them is the $12$-dimensional Zariski closure $H_0$
of the space of twisted cubic curves in $H$. The other $15$-dimensional
component $H_1$ contains the planar cubics with an additional point in $\P_3$.
Piene and Schlessinger showed in \cite{P/S} that $H_0\cap H_1$ is
smooth, rational of dimension $11$ and that the two components
intersect transversally. The curves in the intersection are planar,
singular cubics with an embedded point at one of the
singularities. 

Now we determine the types of the sheaves in $M$ and their resolutions.

We recall first a few facts about locally Cohen-Macaulay curves 
$C\subset\P_3$ of degree $d$:
  \begin{itemize} \label{Facts}
    \item[(C1)] The arithmetic genus satisfies
  $p_a(C)\leq \binom{d-1}{2}$ with equality if and only if $C$ is contained
in a plane $H\subset\P_3$.
    \item[(C2)] If $C$ is not planar then
      $h^1(\ki_C(n))\leq\binom{d-2}{2}-p_a(C)$ for all $n\in\Z$.
    \item[(C3)] $\reg(\ki_C)\leq \binom{d}{2}+1-p_a(C)$, where $\reg$
denotes the Mumford-Castelnuovo regularity.
  \end{itemize}
\vskip5mm

\begin{lemma}\label{custlem}
 For any cubic curve $C$ in $H_0\setminus H_1$ the structure sheaf $\ko_C$
is stable.
\end{lemma}

\begin{proof} The curves in $H_0\setminus H_1$ are precisely the connected
locally Cohen-Macaulay curves of degree $3$ in $\P_3$ which are not contained in a plane
and have Hilbert polynomial $P_{\ko_C}(m)=3m+1$, see e.g. \cite{P/S} or 
\cite{Har}. Given a proper quotient $\ko_C\to \kq$ of pure dimension $1$,
$\kq$ is the structure sheaf $\ko_{C^\prime}$ of a locally Cohen-Macaulay curve 
$C^\prime\subset C$ of $\deg(C^\prime)\le 2$.
If $\deg(C^\prime)=1$, then $C^\prime$ is a line with $R_{\ko_{C^\prime}}(m)=m+1$
and thus $R_{\ko_C} <R_\kq$. If $\deg(C^\prime)=2$, then $C^\prime$ is a 
union of two skew lines or a conic and thus $P_{\ko_{C^\prime}}(m)=2m+2$ or
$2m+1$. In both cases $R_{\ko_C}<R_\kq$.
\end{proof}
\vskip5mm

\begin{lemma}\label{ExtLem}
Let $\kf$ be a stable sheaf on $\P_3$ with Hilbert polynomial $3m+1$. Then
the support $C=Z(\kf)$ of $\kf$ is a connected locally Cohen-Macaulay curve of 
degree 3.     

\begin{itemize}\itemsep2pt
      \item[(i)] If $C$ is contained in a plane
        $H\subset\P_3$ then there exists a non-split extension
     $$0\to\ko_C\stackrel{s}{\longrightarrow}\kf
            \longrightarrow k_p\to 0,$$
     where $p\in C$ and $k_p$ denotes the skyscraper sheaf supported
     on that point.
      \item[(ii)] If $C$ is not contained in a plane, then 
$C\in H_0\setminus H_1\subset\text{\rm Hilb}_{\,3m+1}(\P_3)$ and
        $\kf$ is isomorphic to the structure sheaf $\ko_C$.
        
    \end{itemize}
\end{lemma}

\begin{proof}
Because $P_\kf(0)=1$, the sheaf $\kf$ has a non-zero section $s$. Let $C$ denote
the support of $s$. The corresponding homomorphism $\ko\to \kf$ induces an 
embedding $\ko_C\subset \kf$. Since $\kf$ is purely $1$-dimensional
Cohen-Macaulay sheaf, also $C$ is purely $1$-dimensional
and a locally Cohen-Macaulay curve. We have 
$P_{\ko_C}(m)\le P_{\kf}(m)$ 
for large $m$ and thus $1\le \deg(C)\le 3$. If $\deg(C)\le 2$, 
then $P_{\ko_C}(m)=m+1$, $2m+1$ or $2m+2$, depending on whether $C$ is a line,
a conic or a union of two skew lines. But this is excluded by the stability 
of $\kf$. Therefore, $C$ is a cubic curve. It is connected, 
because any proper connected component of it would violate stability, too.

(i) If $C$ is contained in a plane, then $P_{\ko_C}(m)=3m$.  It follows that 
the cokernel $\kf/\ko_C$ is a skyscraper sheaf $k_p$. 
Because $\kf$ is stable, it cannot be isomorphic to the direct sum 
$\ko_C\oplus k_p.$ This proves (i). Furthermore, the ideal sheaf $\ki_C$ also 
annihilates $\kf$, because, if g is germ in $\ki_C$, then $g\kf$ has 0-dimensional 
support and must vanish because $\kf$ has no 0-dimensional torsion. This proves
that $C=Z(\kf)$ in the planar case.

(ii) Assume now that $Z(\kf)$ is not contained in a plane. Then also
$C$ is not contained in a plane by the previous part (i). Then
$1-P_{\ko_C}(0)=p_a(C)<1$ or $0<P_{\ko_C}(0)$, see (C1). On 
the other hand, the stability of $\kf$ implies
$P_{\ko_C}(0)/3\le 1/3$ and hence $P_{\ko_C}(0)=1$. Therefore
$P_{\ko_C}(m)=3m+1=P_{\kf}(m)$, which yields
$\ko_C=\kf$. Moreover, $C=Z(\kf)$ must be 
contained in $H_0\setminus H_1$. This 
completes the proof. Note, that in both cases $h^0(\kf)=1$.
\end{proof}
\vskip5mm

\begin{lemma}\label{plresl}
Let $\kf$ be a stable sheaf on $\P_3$ with Hilbert polynomial
$3m+1$ and let $C=Z(\kf)$ be its supporting cubic curve. Then the following are 
equivalent:
\begin{itemize}
\item[(i)] There exists a non-split extension of type
  $0\to\ko_C\longrightarrow\kf\longrightarrow k_p\to 0$
  with $p\in C$.
\item[(ii)] $\kf$ is supported on a plane $H$ and has a free resolution
\begin{equation}0\to 2\ko_H(-2)\underset{\binom{q_1\ l_1}{q_2\
    l_2}}{\longrightarrow}\ko_H\oplus\ko_H(-1)\to \kf\to 0,\label{plres}
\end{equation}
\end{itemize}
\end{lemma}

{\bf Remark:} It can easily be verified that the resolution (\ref{plres}) is the
Beilinson resolution of $\kf$ on the plane $H$.

\begin{proof}
(i) to (ii):  The extension sequence implies that $P_{\ko_C}(m)=3m$ and 
hence $p_a(C)=1$. Then $C$ is contained in a plane $H$. Because $C$ is the 
support of $\kf$, the extension sequence is a sequence of $\ko_H$-modules.
By standard homological algebra, a resolution of $\kf$ can be obtained from 
the Koszul resolution of $k_p$ twisted by $\ko_H(-1)$ and the resolution of
$\ko_C$ by its cubic form $f$. The result is the resolution
$$0\to\ko_H(-3)\xrightarrow{B}\ko_H(-3)\oplus 2\ko_H(-2)\xrightarrow{A}
\ko_H \oplus \ko_H(-1)\to\kf\to 0$$
with $$A={\text{\tiny$\left( \begin{array}{cc}
     f & 0\\
     q_1 & l_1\\
     q_2 & l_2
     \end{array}\right)$}}\quad\text{ and }\quad B=(\lambda,-l_2, l_1),$$
such that in particular $\lambda f=q_1\,l_2-q_2\,l_1$.
Suppose $\lambda=0$. Then $l_1 | q_1$
and $l_2 | q_2$ since $l_1$ and $l_2$ are independent. Hence there
exists a linear form $w$ such that $q_1=w\,l_1$ and $q_2=w\,l_2$, and it would 
follow that $A$ is equivalent to
{\tiny $\left(\begin{array}{cc} f & 0\\ 0 & l_1\\0 &
      l_2\end{array}\right)$}.
But then $\kf$ decomposes into
$\kf\cong\ko_C\oplus k_p$ contradicting the stability of
$\kf$. Consequently, we can
assume $\lambda=1$. Then
$$A={\text{\tiny $\left(\begin{array}{cc} q_1\,l_2-q_2\,l_1 & 0\\ q_1 & l_1\\q_2 & l_2
    \end{array}\right)$}}\quad\text{ is equivalent to }\quad
{\text{\tiny $\left(\begin{array}{cc} 0 & 0\\ q_1 & l_1\\q_2 & l_2
    \end{array}\right)$}},$$
and thus one can delete the ghost summand $\ko_H(-3)$.
The homomorphism of the matrix 
{\tiny $\left(\begin{array}{cc} q_1 & l_1\\q_2 & l_2
    \end{array}\right)$} is injective because $q_1\, l_2-q_2\,l_1\neq 0$.

(ii) to (i): Since the cokernel of the composition of the projection
$\ko_H \oplus \ko_H(-1)\to \ko_H(-1)$
with the matrix of the resolution (\ref{plres}) is $k_p$, we obtain
a surjective map
$\kf\to k_p$. It follows easily that its kernel is $\ko_C$.

\end{proof}
\vskip5mm

\begin{lemma}\label{hireslem} Let $C$ be a curve in $H_0\smallsetminus
    H_1\subset\text{\rm Hilb}_{\,3m+1}(\P_3)$. Then $\ko_C$ has a free
    resolution
\begin{equation}\begin{CD}
      0 @>>> 2 \ko_{\P_3}(-3) @>>> 3\ko_{\P_3}(-2)
        @>>> \ko_{\P_3} @>>>\ko_C @>>> 0
      \end{CD}.\label{spres}\end{equation}
\end{lemma}

\begin{proof} Since $C$ is not contained in a plane we get $h^0\ki_C(1)=0$.
The ideal sheaf $\ki_C$ is $4$-regular and $h^1\ki_C(2)$ vanishes,
   cf.~remarks (C1) -- (C3). Hence $P_{\ki_C}(m)=\frac{1}{6}m^3+m^2-
   \frac{7}{6}m$ is equal to $h^0\ki_C(m)$ for all $m\geq 2$.
   This implies that the saturated
   homogeneous ideal $I=H^0_{\ast}(\ki_C)$ of the curve is
   minimally generated by three independent quadratic forms
   $I=(q_1,q_2,q_3)$. The structure sheaf
   $\ko_C$ is stable according to lemma \ref{custlem} and consequently $C$ is
   locally Cohen-Macaulay. Therefore,
   we get an exact sequence
   $$\begin{CD} 0 @>>> \ke @>>> 3\ko_{\P_3} @>>{\text{\tiny
      $\phi = \begin{pmatrix} q_1 \\ q_2 \\ q_3 \end{pmatrix}$}}> \ki_C(2)
       @>>> 0,
   \end{CD}$$
where $\ke$ is a locally free sheaf of rank $2$ on $\P_3$.
$H^1_{\ast}(\ke)=0$ because the map $\phi$ is already surjective on global
sections. Furthermore, (C2) implies $H^2\ke(m)\cong
H^1\ki_C(2+m)=0$ for all
$m\in\Z$. Thus, due to Horrock's criterion, see \cite{OSS}, the vanishing of the
intermediate cohomology implies that the bundle $\ke$ splits,
i.e.~$\ke=\ko_{\P_3}(a)\oplus\ko_{\P_3}(b)$. We know that
$a,b<0$. Using Riemann-Roch, we get $c_1(\ki_C(2))=2$
 since $\text{rk}(\ki_C(2))=1$ and
 $P_{\ki_C(2)}(m)=\frac{1}{6}m^3+2m^2+\frac{29}{6}m+3$.
 It follows that $a=b=-1$ because $0=c_1(\ke)+c_1(\ki_C(2))=a+b+2$.
\end{proof}
\vskip5mm

\begin{proposition}\label{resprop}
   Every stable sheaf $\kf$ on $\P_3$ with $P_{\kf}(m)=3m+1$ has
   a free resolution
$$0\to2\,\ko_{\P_3}(-3)\xrightarrow{B}\ko_{\P_3}(-1)\oplus3\,\ko_{\P_3}(-2)
\xrightarrow{A}\ko_{\P_3}\oplus\ko_{\P_3}(-1)\to\kf\to 0.\label{mainres}\eqno (CR)$$
If $C=Z(\kf)$ is contained in a plane $H=Z(w)$ the matrices can be chosen
     in normal form
     $$B=\begin{pmatrix} -q_1 & -l_1 & w & 0\\ -q_2 & -l_2 & 0 & w\end{pmatrix}\;\text{ and }\;
       A=\begin{pmatrix} w & 0 \\ 0 & w \\ q_1 & l_1 \\ q_2 & l_2
       \end{pmatrix},\,
     \text{where }I_C=(q_1\,l_2-q_2\,l_1,\,w)\text{ and }l_1\sw l_2\neq 0.$$

     If $Z(\kf)$ is not contained in any plane then the normal forms
     are
     $$B=\begin{pmatrix} 0 & l_1 & l_2 & l_3\\ 0 & l_4 & l_5 & l_6
         \end{pmatrix}\;\text{ and }\;
       A=\begin{pmatrix} 0 & 1 \\ q_1 & 0 \\ q_2 & 0 \\ q_3 & 0
       \end{pmatrix}$$
 with linear forms $l_i$ and quadratic forms $q_i.$
\end{proposition}
 \begin{proof} Say, $Z(\kf)\subset H$. Then by lemma \ref{ExtLem}.(i)
    and lemma \ref{plresl}.(ii), there exists a resolution
    \[
      0\to 2\ko_H(-2)\underset{\binom{q_1\ l_1}{q_2\
      l_2}}{\longrightarrow}\ko_H\oplus\ko_H(-1)\to \kf\to 0
    \]
    of $\kf$. Since the equation $w$ of the plane $H$ is annihilated by
    $\kf$ also the map $\ko_{\P_3}(-1)\oplus3\,\ko_{\P_3}(-2)
        \stackrel{A}{\longrightarrow}\ko_{\P_3}\oplus\ko_{\P_3}(-1)$
        has $\kf$ as cokernel. It is then easy to verify that the matrix of relations
    for $A$ is indeed $B$. The non-planar case
follows immediately from lemma \ref{ExtLem}.(ii) and lemma \ref{hireslem}.
  \end{proof}
\vskip5mm

{\bf Remark:} From the resolution (CR) in \ref{resprop} we see that the
Castelnuovo-Mumford regularity of any stable sheaf with
Hilbert polynomial $3m+1$ is equal to $1$ and hence $H^1\kf=0$.
This fact and the Beilinson-II spectral sequence
$$H^s(\P_3,\kf\otimes\Omega^{-r}_{\P_3}(-r))\otimes_k
\ko_{\,\P_3}(r) =: E_1^{rs} \Longrightarrow E_\infty^i\cong\left\{
{\text{gr}(\kf),\atop 0,}{\text{ for }
    i=0\atop\text{ otherwise }}\right.$$
induce an exact sequence
\[0 \to H^1\kf(-1) \otimes \ko_{\P_3}(-3)
  \to H^0(\kf\otimes\Omega_{\P_3}^1(1)) \otimes \ko_{\P_3}(-1)
  \oplus H^0(\kf\otimes\Omega_{\P_3}^2(2))\otimes\ko_{\P_3}(-2)\]
\begin{equation}\to H^0\kf \otimes \ko_{\P_3} \oplus H^1(\kf\otimes\Omega_{\P_3}^1(1))
  \otimes \ko_{\P_3}(-1) \to \kf \to 0\quad\label{beilinson}\end{equation}
because the cohomology groups $H^0\kf(-j),\; j>0$, vanish due to the
stability of $\kf$.

Using the two sequences (\ref{plres}) and (\ref{spres}), it is an easy exercise
in cohomology to check that
\[
  \begin{array}{lcl}
    h^0(\kf\otimes\Omega^3(3))=0 &, & h^1(\kf\otimes \Omega^3(3))=2 \\
    h^0(\kf\otimes\Omega^2(2))=0 & , & h^1(\kf\otimes \Omega^2(2))=3\\
\multicolumn{3}{l}
{    h^0(\kf\otimes \Omega^1(1))=h^1(\kf\otimes \Omega^1(1))=
          \left\{ \,
          \begin{matrix} 1 & \text{ if }Z(\kf)\text{ is planar.} \\
                         0 & \text{otherwise.}
          \end{matrix}\right.}
  \end{array}
\]
Hence, the resolution (CR)  is the Beilinson
resolution (\ref{beilinson}) if the sheaf is supported on a plane. In
the case of non-planar sheaves, the terms
$\ko_{\P_3}(-1)$ can be cancelled from (CR) to
give the Beilinson resolution (\ref{beilinson}) in that case.
Furthermore, due to the explicit construction in this section
we also know {\sl exactly which types of matrices $A$ and $B$ can occur}.

The relative version of Beilinson's theorem, see \cite{OSS}, p.~306, specializes
in our case to
\vskip5mm

\begin{proposition}\label{relbeil}
1) Let $\bkf\in\coh(S \times \P_3 )$ be an $S$-flat family of semi-stable sheaves
$\bkf_s$ with Hilbert polynomial $P_{\bkf_s}(m)=3\,m+1$ for all $s\in
S$. Then there is an exact sequence

$0\to\kb_2\boxtimes\ko_{\P_3}(-3)\to\ka_1\boxtimes\ko_{\P_3}(-1)\oplus\kb_1\boxtimes
\ko_{\P_3}(-2)$

\hfill $\to \ka_0\boxtimes\ko_{\P_3}\oplus\kb_0\boxtimes\ko_{\P_3}(-1)\to\bkf\to 0$

with
\[
\begin{array}{ll}
                     & \kb_2=R^1{p_\ast}(\bkf\otimes\Omega^3(3))\\[1.0ex]
\ka_1={p_\ast}(\bkf\otimes \Omega^1(1)) &
\kb_1=R^1{p_\ast}(\bkf\otimes\Omega^2(2))\\[1.0ex]
\ka_0={p_\ast}\bkf & \kb_0=R^1{p_\ast}(\bkf\otimes \Omega^1(1))
\end{array}
\]

2) The sheaves $\kb_2, \kb_1, \ka_0$ are locally free of rank $2,3,1$
respectively.

3) The sheaves $\ka_1$ and $\kb_0$ are supported on the subscheme $S_1 \subset
S$ of points $s\in S$ for which $\bkf_s$ is planar.
\end{proposition}
\vskip5mm
\begin{proof}
The first part follows from the relative Beilinson spectral sequence and the vanishing
of cohomology groups $H^0(\bkf_s(-j))$ of the stable sheaves $\bkf_s,
s\in S$. The
second part follows from the base change theorem for the points of $S$ using
that $h^1(\bkf_s\otimes \Omega^3(3))=2$, $h^1(\bkf_s\otimes \Omega^2(2))=3$ and
$h^0(\bkf_s)=1$ for any $s$ and the vanishing $h^0(\bkf_s\otimes
\Omega^3(3))=h^0(\bkf_s\otimes \Omega^2(2))=0$ (the scheme $S$ might be
non--reduced).

Because $\bkf$ is $S$--flat, the function $s\,\mapsto\,h^0(\bkf_s\otimes
\Omega^1(1))=h^1(\bkf_s\otimes \Omega^1(1))$ is upper semi-continuous and
defines a reduced, closed subscheme $S_1\subset S$ as its support. Then
$S_1=Z(\ka_1)=Z(\kb_0)$.
\end{proof}

There arise two non--standard problems from the last two propositions: The
resolution of type (CR)
is not the Beilinson resolution if the sheaf is the structure sheaf $\ko_C$
of a non-planar cubic.
Therefore, the Beilinson construction cannot be used directly in order to
describe the parameters for the moduli space. Secondly, the
automorphism group of a complex (CR) is
not reductive. Both problems are overcome in section 5.
\vskip10mm

\section{The two components}

\rm Let $M=M_{3m+1}(\P_3)$. As stated in proposition \ref{genmod}, there is a
universal sheaf $\bku$ on $M\times \P_3$. By proposition \ref{relbeil}, 3),
there exists
a closed subscheme $M_1\subset M$ with reduced structure  consisting
of all sheaves with planar support. We will show in section 6 that
$M_1$ is irreducible, smooth and rational
of dimension $13$. Let now $M_0$ denote the closure of the open subscheme
$M\smallsetminus M_1$ of non--planar sheaves.  Then we have
\[
M=M_0\cup M_1.
\]
\vskip5mm

\begin{lemma}\label{spcub} $M\smallsetminus M_1$ is isomorphic to
$H_0\smallsetminus H_1$, the open
  part of non--planar cubics of the Hilbert scheme $\hilb_{3m+1}(\P_3)$.
\end{lemma}
\vskip5mm

\begin{proof}
There is a set-theoretical bijection $\phi\,:\, M_0\smallsetminus M_1
  \longrightarrow H_0\smallsetminus H_1$ due to lemma \ref{ExtLem}, (ii). This
  bijection is indeed an isomorphism:
Let $Z_0\subset (H_0\smallsetminus H_1)\times \P_3$
be the universal cubic of $\hilb_{3m+1}(\P_3)$ restricted to
  $H_0\smallsetminus H_1$. Then there is a resolution
\[
0\to \ke_1\boxtimes\ko_{\P_3}(-3)\to \ke_0
\boxtimes \ko_{\P_3}(-2)\to \ko_{(H_0\smallsetminus H_1)\times\P_3}
{\longrightarrow} \ko_{Z_0} \to 0
\]
with $\ke_1$ and $\ke_0$ locally free on $H_1\smallsetminus H_0$ of
rank $2$ and $3$. It restricts to the resolution from
lemma \ref{hireslem} along the fibers.
On the other hand, let $\bku_0$ denote the restriction of the
universal sheaf $\bku$ to
  $(M\smallsetminus M_1)\times \P_3$. Since the sheaves $\ka_1$ and
$\kb_0$ in the resolution of $\bku$ vanish on $M\smallsetminus M_1$
  (cf. proposition \ref{relbeil}), we get
\[
0\to (\kb_2\otimes \ka_0^\ast)\boxtimes\ko_{\P_3}(-3)\to (\kb_1\otimes
\ka_0^\ast )\boxtimes \ko_{\P_3}(-2)\to \ko_{(M\smallsetminus M_1)\times\P_3}
{\longrightarrow} \bku_0\otimes p^\ast\ka_0^\ast\to 0.
\]

Because $\bku_0\otimes p^\ast\ka_0^\ast$ represents a family of cubic curves in 
$H_0\smallsetminus H_1$, it follows from the universal 
property of H and its open subset $H_0\smallsetminus H_1$, 
that there is a unique morphism
$M_0\smallsetminus M_1\xrightarrow{\alpha}H_0\smallsetminus H_1$ such that
$\bku_0\otimes p^\ast\ka_0^\ast$ is the pullback of $\ko_{Z_0}$. Since on the 
other hand $\ko_{Z_0}$ is also a family of stable sheaves belonging to 
$M_0\smallsetminus M_1$, there is a unique morphism 
$H_0\smallsetminus H_1\xrightarrow{\beta}M_0\smallsetminus M_1$
such that conversely $\ko_{Z_0}$ is the pullback of 
$\bku_0\otimes p^\ast\ka_0^\ast$. By uniqueness $\alpha$ and $\beta$ are 
inverse to each other. The underlying map of $\alpha$ is exactly the 
bijection $\phi$ above.
\end{proof}

Piene and Schlessinger's result \cite{P/S} on
$H_0\subset\hilb_{3m+1}(\P_3)$ implies
then the following corollary.
\begin{corollary}\label{spcub1} The subscheme $M_0$ is an irreducible
  component of $M$ of dimension $12$ and smooth along $M_0\smallsetminus M_1$.
\end{corollary}

We will show in section 8 that the component $M_0$ as a whole is
isomorphic to $H_0$ and that it can be identified with
a blowup of the space of nets of quadrics
as described in \cite{E/S} such that $M_0\cap M_1$ is the exceptional
divisor.
\vskip10mm

\section{A parameter space}
The Beilinson resolution of proposition \ref{relbeil} over $M\times\P_3$ of
the universal sheaf contains the sheaves $\ka_1$ and $\kb_0$ which are
supported on $M_1$ and so are not locally free on $M$. On the other hand, they
do not cancel because the homomorphism between them is zero.
In the following, we construct a quasi--affine variety $X$
consisting of the matrices which occur in the resolution in proposition
\ref{resprop} such that the moduli space $M$ is a geometric quotient of $X$
and such that the variety $X$ is the minimal one allowing a locally free
resolution of the lifted universal sheaf. To begin with, we prove the
following
\vskip5mm

\begin{lemma}\label{locpres} Any $\kf$ in $M$ has an open
  neighbourhood $U\subset M$ such that there is a resolution
\[
0\to 2\widetilde{\ko}_U(-3)\to \widetilde{\ko}_U(-1)\oplus
3\widetilde{\ko}_U(-2)\to\widetilde{\ko}_U\oplus\widetilde{\ko}_U(-1)\to
\bku |_{U\times\P_3}\to 0,
\]
where $\widetilde{\ko}_U(d)=\ko_U\boxtimes\ko_{\P_3}(d)$ and $\bku$
denotes the universal sheaf on $M\times\P_3$.
\end{lemma}

\begin{proof} Consider the relative Beilinson resolution of $\bku$
from
proposition \ref{relbeil}. If $\kf\in M_0\smallsetminus M_1$ choose
an open neighbourhood $U\subset M_0\smallsetminus M_1$ of $\kf$ and
add the complex $0\to \widetilde{\ko}_U(-1)=\widetilde{\ko}_U(-1)\to
0$ to the restriction of the resolution to $U\times\P_3$. In the
case $\kf\in M_1\smallsetminus M_0$, the restriction of the
resolution to some open neighbourhood $\kf\in U\subset M_1
\smallsetminus M_0$ is already of the proposed type since
$\ka_1$ and $\kb_0$ are locally free on their support $M_1$.

So let $\kf \in M_0\cap M_1$. There exists an open neighbourhood
$U_1\subset M_1$ of $\kf$ with
\[0\to 2\tko_{U_1}(-3)\xrightarrow{\widetilde{B}_1}
\tko_{U_1}(-1)\oplus
3\tko_{U_1}(-2)\xrightarrow{\widetilde{A}_1}\tko_{U_1}\oplus\tko_{U_1}(-1)\xrightarrow{\widetilde{\varphi}_1}\bku|_{U_1\times\P_3}\to
  0\ .
\]
Here the component of $\widetilde{A}_1$ from $\tko_{U_1}(-1)$ to
$\tko_{U_1}(-1)$ is zero. We are going to extend this resolution to a
resolution of $\bku$ on an open set $U$ in $M$ with
$U\cap M_1=U_1$.

For that let $\ki$ denote the ideal
sheaf of $M_1\subset M$. The exact sequence $0\to p^\ast\ki \otimes
\bku\to\bku\to\bku|_{M_1\times\P_3}\to 0$ induces
$$ H^0(U\times\P_3,\bku(i)) \longrightarrow
H^0(U_1\times\P_3,\bku(i)|_{M_1\times\P_3})
\longrightarrow H^1(U\times\P_3,p^\ast\ki\otimes\bku(i)),\quad i=0,1$$
where $U\subset M$ is some open affine with $U\cap
( M_0\smallsetminus M_1 )\neq\emptyset$ and $U\cap M_1=U_1$.

{\bf Claim}: For any coherent sheaf $\kg$ on $M$, $R^1 p_{\ast}(p^{\ast}
\kg\otimes\bku (i))=0$ for $i=0,1$.

Note that $R^1 p_\ast\bku(i)=0$ since $\reg(\kf)=1$ for all $\kf\in M$.
If $\kg$ is locally free the claim follows from the projection
formula: $R^1p_\ast(p^\ast\kg\otimes\bku(i))=\kg\otimes
R^1p_\ast\bku(i)$. Otherwise take a free resolution
$\ke_\bullet\to\kg\to 0$ and split it into short exact sequences. Then
apply $p^\ast$ to the building blocks and tensor them with
$\bku(i)$. Due to the flatness, exactness is preserved. A look
at the long exact sequence associated to $p_\ast$ finishes the proof
of the claim.

In particular, we obtain
$$H^1(U\times\P_3,p^\ast\ki\otimes\bku(i))\cong H^0(U,R^1p_\ast(p^{\ast}
\ki\otimes\bku(i))=0$$
from Leray's spectral sequence and the claim above. It follows that
$\widetilde{\varphi}_1$ can be extended:
\[
\xymatrix{
0\ar[r] & \kk_1 \ar[r] & \tko_{U_1}\oplus\tko_{U_1}(-1) \ar[r]^-
{\widetilde{\varphi}_1} & \bku|_{U_1\times\P_3} \ar[r] & 0 \\
 & &
\tko_{U}\oplus\tko_{U}(-1) \ar[r]^-{\widetilde{\varphi}} \ar@{->>}[u]^-
{\text{res}}\ar[r] &
\bku|_{U\times\P_3}\ar@{..>}[r] \ar@{->>}[u] & 0.}
\]
After possibly shrinking $U$, we can assume that
$\widetilde{\varphi}$ is surjective. Let $\kk$ denote the kernel of
$\widetilde{\varphi}$. By flatness,
$\kk\otimes
\ko_{M_1}=\kk_1=\im(\text{res}|_\kk)$. The process of extending
$\widetilde{\varphi}_1$ can be repeated for the situation
\[
\tko_{U_1}(-1)\oplus
3\tko_{U_1}(-2)\xrightarrow{\widetilde{A}_1} \kk_1\to 0
\]
in order to obtain the first part of the resolution in the lemma. The last
part of the resolution is obtained by the same procedure applied to the kernel
of $\tko_{U_1}(-1)\oplus 3\tko_{U_1}(-2)\to {\kk_1}$.
\end{proof}
\vskip5mm

\begin{sub}\label{parasp}\rm{\bf The Parameter space.}\\
The space of parameters for $M$ will be the space of resolutions described in
proposition \ref{resprop}. To be precise, let $\P=\P V=\P_3$ and let
\[
W\subset\hom(2\ko_{\P}(-3), \ko_{\P}(-1)\oplus 3\ko_{\P}(-2))\times
\hom(\ko_{\P}(-1)\oplus
3\ko_{\P}(-2), \ko_{\P}\oplus \ko_{\P}(-1))
\]
be the locally closed subvariety of pairs $(B, A)$ for which the induced
sequence
$$
0\to k^2\xrightarrow{\widetilde{B}} S^2V^\ast\oplus(k^3\otimes
V^\ast)\xrightarrow{\widetilde{A}} S^3V^\ast\oplus S^2V^\ast,\quad
\widetilde{B}=H^0B(3),\;\widetilde{A}=H^0A(3)
\eqno(E)
$$
is exact. $W$ is acted on by the automorphism group
\[
G=\aut(2\ko_\P(-3))\times \aut(\ko_\P(-1)\oplus 3\ko_\P(-2))\times
\aut(\ko_\P\oplus \ko_\P(-1)),
\]
the non--reductive group of triples of matrices
\[
g_1,\quad g_2=\left(
  \begin{array}{c|ccc}
\alpha & 0 & 0 & 0\\\hline
u_1 & & &\\
u_2 & & g &\\
u_3 & & &
  \end{array}\right),\quad g_3=\left(
  \begin{array}{cc}
\beta & 0\\
u & \gamma
  \end{array}\right)
\]
with $g_1\in \GL_2(k)$, $g\in \GL_3(k)$, $u, u_i\in V^\ast$ and
$\alpha,\beta,\gamma\in k^\ast$, the action being given by
$(g_1 B g_2^{-1}, g_2 Ag_3^{-1}).$
The pair $(B, A)\in W$ is called {\bf stable} if
$$A=\begin{pmatrix} z & \lambda \\ q_1 & z_1 \\ q_2 & z_2 \\ q_3 & z_3
       \end{pmatrix}$$
satisfies $\lambda\neq 0$ or $z_1\sw z_2\sw z_3\neq 0$. We let
$X\subset W$ be the open subset of stable pairs.
\end{sub}
\vskip5mm

\begin{lemma}\label{lempara} For any $(B,A)\in X$ the sheaf $\kf_A=\Coker(A)$
is stable and has the resolution by $(B, A)$ as in proposition \ref{resprop}.
\end{lemma}

\begin{proof} The pair $(B, A)$ defines the complex
$$
0\to 2\ko_\P(-3)\xrightarrow{B}\ko_\P(-1)\oplus
3\ko_\P(-2)\xrightarrow{A}\ko_\P\oplus\ko_\P(-1) \to\kf_A\to 0\eqno(S)
$$
which is already exact except possibly at $\ko_\P(-1)\oplus 3\ko_\P(-2)$, because
$\widetilde{B}$ has rank $2$. In order to prove exactness and stability, we
distinguish the cases $\lambda\neq 0$ and $\lambda=0$ for the degree $0$ entry
of $A$. If $\lambda\neq 0$, it follows that $(B, A)$ is $G$--equivalent to a
pair
\[
\begin{array}{ll}
\left(  \begin{array}{c|lll}
0 & l_1 & l_2 & l_3\\
0 & l_4 & l_5 & l_6
  \end{array}\right), & \left(
  \begin{array}{cc}
0 & 1\\
q_1 & 0\\
q_2 & 0\\
q_3 & 0
  \end{array}\right)
\end{array}.
\]
In that case the matrix $B$ of rank 2 cannot be equivalent to one with
$l_3=l_6=0$, because then $q_1=q_2=0$ and
$\widetilde{A}$ would have a higher dimensional kernel. Then the linear part
$B'$ of $B$ is a stable point of $\hom(k^2, k^3\otimes V^\ast)$ for the action
of $\GL_2(k)\otimes \GL_3(k)$, see \cite{E/S} and \ref{nets}. It follows that
the forms $q_i$
are the minors of this matrix $B'$, that $\kf_A\cong \ko_C\cong\ko_\P/(q_1,
q_2, q_3)\ko_\P$ is stable and that $(S)$ is exact. If $\lambda=0$, one proves by
an elementary procedure that $(B,A)$ is $G$--equivalent to a pair
\[
\left(
  \begin{array}{llcc}
-q_1 & -l_1 & w & 0\\
-q_2 & -l_2 & 0 & w\\
  \end{array}\right),\quad\left(
  \begin{array}{rr}
w & 0\\
0 & w\\
q_1 & l_1\\
q_2 & l_2
  \end{array}\right)
\]
with $w\sw l_1\sw l_2\neq 0$. Then also $q_1l_2-q_2l_1\neq 0,$ because
otherwise the sequence (E) would not be exact. In this case
the sheaf $\kf_A$ has support contained in the plane $H=Z(w)$ and has
the resolution
$$0\to 2\ko_H(-2)\underset{\begin{pmatrix} q_1 & l_1\\q_2 &
    l_2\end{pmatrix}}{\longrightarrow} \ko_H\oplus \ko_H(-1)\to\kf_A\to 0. $$
Therefore $\kf_A$ is one of the stable sheaves of the moduli space
$M_{3m+1}(H),$
see \ref{plcase} and \ref{gqu1}. Moreover, due to the shape of the
matrix pair, the complex $(S)$ is also exact in this case.
\end{proof}
\vskip5mm

The lemma above shows that
$ X\longrightarrow M,\; (B,A)\mapsto[\kf_A]$
is $G$--equivariant and
surjective. It follows from the
existence of the universal family constructed next
that the map is also a morphism.
\vskip5mm

\begin{sub}\label{univsh}\rm{\bf The universal sheaf on $X$.}\\
On $H_1:=\hom(2\ko_\P(-3), \ko_\P(-1)\oplus 3\ko_\P(-2))$, we are given
the tautological homomorphism
\[
k^2\otimes \ko_{H_1}\xrightarrow{b}(S^2V^\ast\oplus k^3\otimes V^\ast)\otimes
\ko_{H_1}
\]
and on $H_2:=\hom(\ko_\P(-1)\oplus
3\ko_\P(-2), \ko_\P\oplus \ko_\P(-1))$ we have the tautological maps
\[
\begin{array}{lcll}
\ko_{H_2} & \xrightarrow{a_1} & (V^\ast\oplus k)\otimes \ko_{H_2} & \text{(first component)}\\
k^3\otimes \ko_{H_2} & \xrightarrow{a_2} & (S^2V^\ast\oplus V^\ast)\otimes
  \ko_{H_{2}} & \text{(second component)}\ .
\end{array}
\]
Combining these three maps with the corresponding evaluation homomorphisms
$S^mV^\ast\otimes \ko_\P\to
\ko_\P(m)$, we obtain homomorphisms
\[
k^2\otimes \ko_{H_1}\boxtimes \ko_\P(-3)\xrightarrow{\widetilde{B}}
\ko_{H_1}\boxtimes (\ko_\P(-1)\oplus k^3\otimes \ko_\P(-2))
\]
and
\[
\ko_{H_2}\boxtimes (\ko_\P(-1)\oplus k^3\otimes
\ko_\P(-2))\xrightarrow{\widetilde{A}}\ko_{H_2}\boxtimes(\ko_\P\oplus \ko_\P(-1)).
\]
Finally, restricting to $X\times \P\subset W\times \P\subset H_1\times
H_2\times \P$ and using the notation
$\tko_X(m)=\ko_X\boxtimes \ko_\P(m)$, we get the complex
$$
0\to 2\tko_X(-3)\xrightarrow{\widetilde{B}} \tko_X(-1)\oplus
3\tko_X(-2)\xrightarrow{\widetilde{A}}\tko_X\oplus\tko_X(-1)\to\bkf\to 0\eqno (R)
$$
with $\bkf=\coker(\tilde{A})$. It is exact after restriction to a
point $(B,A)$ of $X$, giving
$\bkf_{(B,A)}=\coker(A)=\kf_A$. This proves that the complex (R) is exact and that
$\bkf$ is a flat family of sheaves over $X$. Now the above family defines a
unique morphism $X\overset{\mu}{\twoheadrightarrow} M$
such that $\bkf\cong\mu_{\,\P}^\ast\,\bku$, whose underlying map is
$(B,A)\mapsto [\kf_A]$, which is surjective. It is easy to verify that the
fibers of $\mu$ coincide with the $G$--orbits.
\end{sub}

{\bf Remark}: The variety $X$ consists of two irreducible components
$X_1=\mu^{-1}(M_1)$ and $X_0=\mu^{-1}(M_0)$ with $X_1=\{(B,A)\in X\ |\
\lambda=0\}$ and $X_0$ the closure of $X\smallsetminus X_1$.
\vskip5mm

\begin{theorem}\label{geomqu} The morphism $X\xrightarrow{\mu} M$ is a geometric quotient of
  $X$ by the (non--reductive) group $G$ in the sense of Geometric Invariant
  Theory, see \cite{Mum}, \cite{D/T}.
\end{theorem}
\begin{proof} Since the fibers of $\mu$ coincide with the $G$--orbits, it is
  sufficient to prove that $\mu$ admits local sections (slices) through any
  point $(B,A)$ in $X$. Let $(B,A)\in X$ be given and let
  $a:=[\kf_A]\in M$. By lemma
  \ref{locpres} there is an open neighbourhood $U(a)\subset M$ with a
  resolution
\[0\to 2\tko_{U}(-3)\stackrel{\Psi}{\longrightarrow}
\tko_U(-1)\oplus3\tko_U(-2)\stackrel{\Phi}{\longrightarrow}\tko_U\oplus
\tko_U(-1)\to\bku|_{U\times \P}\to 0
\]
of the universal sheaf. Define a morphism $s : U \rightarrow X$ by
$b \mapsto (\Psi(b),\Phi(b))$. Since $\bku_a\cong\kf_A$ we have
$\Phi(a)\sim A$. Then there is an element
$g\in G$ with $g\,.\,s(a)=(B,A)$ and $g\,.\,s$ is then a section of $\mu$
through $(B,A)$.
\end{proof}

{\bf Remark 1}: The definition of a geometric quotient in \cite{Mum} does not
include that the quotient map has to be affine. In our case we do not
know whether the
action is proper or whether the quotient map is affine.

{\bf Remark 2}: There is no common stabilizer for all the points of $X$. The
stabilizers for the points of $X_1$ are all the same and isomorphic to
$k^\ast\times k$, whereas the stabilizers for the points of
$X_0\smallsetminus X_1$ are all isomorphic to $k^\ast\times k^\ast.$
Therefore, the action of $G$ modulo a subgroup cannot be free and
$X\xrightarrow{\mu} M$ is not a principal fibration.

{\bf Remark 3}: In \cite{D/T} polarizations $\Lambda$ have been introduced by
which
open sets $W^s(G,\Lambda)\subset W^{ss}(G, \Lambda)\subset W$ of stable and
semi-stable points for a {\sl non--reductive} group $G$ can be defined in such
a way that $W^{ss}(G, \Lambda)$ admits a good and projective quotient
$W^{ss}(G, \Lambda)//G$. Moreover, $W^s(G, \Lambda)$ and closed
invariant subsets of $W^s(G, \Lambda)$ have
a geometric quotient.

In the above case $X\subset W$, one can prove that
there exists no polarization $\Lambda$ in the sense of \cite{D/T},
such that $X$
is a closed subset of $W^s(G, \Lambda)$, cf.~\cite{HGF}.
\vskip10mm

\section{The component $M_1$}

Intuitively, the component $M_1$ of $M_{3m+1}(\P_3)$
should be fibered by the schemes $M_{3m+1}(H)$ of stable sheaves on a plane
$H$ in $\P_3$, cf. proof of lemma \ref{lempara}.
In fact, $M_1$ is isomorphic to the relative
Simpson scheme over $\P_3^\ast$ as we will see.
Therefore, we first digress to the moduli
space $M_{3m+1}(\P_2)$, see also \cite{LePot2}.
\vskip5mm

\begin{sub} \label{plcase}\rm {\bf The moduli space $M_{3m+1}(\P_2)$.}\\
Let $M(\P_2)=\3m (\P_2)$. Because of the multiplicity $3$, this scheme is a fine
moduli space with a universal sheaf $\bkv$ over $M(\P_2)\times \P_2$ and it consists
entirely of stable sheaves. By lemma \ref{plresl}, any sheaf $\kf$ in
$M(\P_2)$ has a resolution
$$
0\to 2\ko_{\P_2}(-2)\underset{\begin{pmatrix} q_1 & l_1\\q_2 &
    l_2\end{pmatrix}}{\longrightarrow} \ko_{\P_2}\oplus \ko_{\P_2}(-1)\to\kf\to
    0 \eqno (PR)
$$
with independent linear forms $l_1, l_2$ and quadratic forms $q_1, q_2$ such
that $f=q_1l_2-q_2l_1\neq 0$. The plane cubic $C=Z(\kf)$ has the
equation $f$.
\end{sub}
\vskip5mm

\begin{lemma}\label{singpts} With the above notation let $p\in C$ be the point
  determined by $l_1(p)=l_2(p)=0$. Then the following are equivalent for the
  stalk $\kf_p$.
  \begin{enumerate}
  \item [(i)]  $\kf_p$ is not free.
  \item [(ii)] $p$ is a singular point of $C$.
  \item [(iii)] $q_1$ and $q_2$ vanish both at $p$.
  \end{enumerate}
\end{lemma}

The elementary proof is left to the reader.

{\bf Remark:} This lemma shows that the sheaf $\kf$ is a
line bundle on $C$, even if the support $C$ is singular, if and only
if the point $p$ is not a singular point. If $p$ is a singular point,
$\kf$ is a Cohen--Macaulay sheaf on $C$ of rank $1$. We call the
sheaves $\kf$ in $M(\P_2)$ that are not locally free on their support
{\bf singular} sheaves.

Now we consider a parameter space $Y$ for
$M(\P_2)$, using the same method as in section $5$. We define $Y$ to be the
quasi--affine variety of $2\times 2$ matrices as in (PR). It is acted on by the
group
\[
G = \aut(2\ko_{\P V}(-2))\times \aut(\ko_{\P V}\oplus\ko_{\P V}(-1))
  = \gl_2(k)\times\left.\left\{\begin{pmatrix} \alpha & z\\0 &
  \beta\end{pmatrix}\;\right|\;\alpha\beta\neq 0, z\in V^\ast\right\}
\]
via $(g,h)\,.\,A=gAh^{-1}.$

It is clear from the description (PR) of the
sheaves $\kf$ in $M(\P_2)$ that the isomorphism classes $[\kf]$ correspond
bijectively to the orbits of this action on $Y$. So we expect that
$Y/G=M(\P_2),$ if the quotient exists. Indeed, as in \ref{univsh} we can
construct a
universal sheaf $\bkg$ over $Y\times\P V$ with a presentation
\[
0\to k^2\otimes \ko_Y\boxtimes \ko_{\P V}(-2)
\xrightarrow{\Phi}\ko_Y\boxtimes\ko_{\P V}
\oplus\ko_Y\boxtimes\ko_{\P V}(-1)\to\bkg\to 0.
\]
For any $y\in Y$, the restricted homomorphism $\Phi_y=A$ is one of the matrices
above and $\bkg_y$ belongs to $M(\P_2)$. Thus $\bkg$ is a flat family and
defines a surjective morphism $ Y\xrightarrow{\nu}M(\P_2)$
with $\bkg\cong \nu_{\P_2}^\ast\bkv$. By the same method as in the proof of
theorem \ref{geomqu} one can directly construct local slices of $\nu$ using the
relative Beilinson resolution of the universal sheaf $\bkv$, which is locally
free in this case. We so obtain the
\vskip5mm

\begin{proposition}\label{gqu1} The scheme $M(\P_2)=\3m(\P_2)$ together with
the morphism $ Y\xrightarrow{\nu} M(\P_2)$ is a geometric quotient of $Y$ by
the action of the (non-reductive) group $G$.
\end{proposition}

Let now $Z\subset \P S^3 V^\ast\times \P V$ denote the universal cubic given by
pairs $(\langle f\rangle, \langle x\rangle)$ with $f(x)=0$. The map
\[
\begin{array}{lcl}
\begin{pmatrix}q_1 & l_1\\q_2 & l_2\end{pmatrix} & \mapsto & (
  \langle q_1 l_2 - q_2 l_1\rangle, \langle
l_1\sw l_2\rangle)
\end{array}
\]
is a surjective morphism $Y\xrightarrow{\gamma}Z,$ where we identify
$\langle l_1\sw l_2\rangle$ with the point
$p$ defined by $l_1(p)=l_1(p)=0$ via $V\cong \Lambda^2V^\ast$.
It is not difficult to see that $\gamma$ is
a geometric quotient, cf. also \cite{D/T}. Together with proposition \ref{gqu1}, this implies
\vskip5mm
\begin{corollary}\label{lgqu1} $M_{3m+1}(\P_2)$ is isomorphic to the
  universal cubic $Z\subset \P S^3 V^\ast\times \P V$.
\end{corollary}
\vskip3mm

{\bf Remark 1}: The isomorphism $\3m(\P_2)\cong Z$ has already been mentioned in
\cite{LePot2}.
\vskip3mm

{\bf Remark 2}: A universal sheaf $\bkf$ on $Z\times \P_2$ can also directly be
defined as follows. Let $H=\P S^3 V^\ast$. The dual of the ideal sheaf
sequence $0\to\ki_{\Delta}\to\ko_{Z\times_H Z} \to\ko_\Delta\to 0$ of the
diagonal $\Delta\subset Z\times_H Z$ is
\[
0\to\ko_{Z\times_H Z}\to\bkf\to\ke xt^1(\ko_\Delta, \ko_{Z\times_H Z} )\to 0.
\]
The $\ke xt$--sheaf is in this case isomorphic to $\ko_\Delta$ and the sequence
does not split. For a single cubic we get back the sequence $0\to\ko_C\to\kf_C\to
k_p\to 0$, see lemma \ref{ExtLem}. This construction is motivated by the classical
construction of the Poincar\'{e} bundle for a single elliptic curve.
Using the embedding $Z\times_H Z\subset Z\times \P_2$ determined by 
$((C,x),(C,y))\mapsto ((C,x),y)$, the sheaf $\bkf$ can be considered as a sheaf 
on $Z\times \P_2$ which is flat over $Z$. It is isomorphic to the universal 
sheaf.
\vskip3mm

{\bf Remark 3}: One has $Y=W^s(G,\Lambda)=W^{ss}
(G, \Lambda)$ for the polarization $\Lambda=(\frac{1}{2};\frac{3}{4},
\frac{1}{4})$ in the sense of \cite{D/T}, cf.~\cite{HGF}.
\vskip10mm

\begin{sub}\label{relc}\rm {\bf The relative case.}\\
We prove now that the component $M_1$ is fibered by the spaces $M_{3m+1}(H),$
$H$ a plane in $\P_3.$
Let $S:=\P_3^\ast=\P V^\ast$ and let $\P\kh\to S$ be the
bundle of planes, induced by the tautological subbundle
$\kh\subset V\otimes\ko_S.$ We denote by
\[
M_S=M_{3m+1}(\P\kh/S)\to S
\]
the relative Simpson scheme, such that each fiber $M_s\cong
M_{3m+1}(\P\kh_s)$. See \cite{H/L}, theorem 4.3.7, for the general situation.
This scheme $M_S$ is again a fine moduli space with a universal
sheaf
\[
\bkv_S\ \text{ over }\  M_S\times_S\P\kh\subset M_S\times \P V.
\]
Since $\P\kh$ is locally trivial over $S$, also $M_S$ is locally trivial over
$S$ with fiber $M_{3m+1}(\P_2)$. Corollary \ref{lgqu1} implies that $M_S$
is smooth, irreducible and rational of dimension $13$.
Consider now the relative universal cubic
\[
Z_S\subset \P(S^3 \kh^\ast)\times_S \P\kh.
\]
It is easy to show that there exists an isomorphism
$M_S\cong Z_S$ which extends the isomorphisms
$M_s\cong Z_s$ of the fibers from \ref{lgqu1}.
\end{sub}
\vskip5mm

\begin{proposition}\label{m1} The subscheme $M_1\subset M$ is isomorphic to
  $M_S$ and $Z_S$ and therefore smooth, rational of dimension $13$.
\end{proposition}

\begin{proof} Let $\bku_1$ be the restriction of the universal sheaf $\bku$ to
  $M_1\times \P V$. It has a presentation
\[
\ka_1\boxtimes \ko_{\P V}(-1)\oplus \kb_1\boxtimes
\ko_{\P V}(-2)\xrightarrow{A}\ka_0\boxtimes \ko_{\P V}\oplus \kb_0\boxtimes
\ko_{\P V}(-1)\to\bku_1\to 0
\]
with locally free sheaves $\ka_1, \kb_1, \ka_0, \kb_0$ on $M_1$ of ranks
$1,3,1,1$ respectively, see proposition \ref{relbeil}. Note that
the component $\ka_1\boxtimes\ko_{\P V}(-1)\to\kb_0\boxtimes\ko_{\P
  V}(-1)$ of the
homomorphism $A$ vanishes identically and that the component
$\ka_1\boxtimes \ko_{\P V}(-1)\to \ka_0\boxtimes\ko_{\P V}$
corresponds to a map
\[
\ko_{M_1}\xrightarrow{\omega}\ka_0\otimes \ka_1^\ast\otimes V^\ast.
\]
This section $\omega$ of $\ka_0\otimes \ka_1^\ast\otimes V^\ast$ is
nowhere zero, because $A$ can be expressed locally by a matrix
\[
\left(
  \begin{array}{ll}
w & 0\\
q_1 & z_1\\
q_2 & z_2\\
q_3 & z_3
  \end{array}
\right)
\]
of linear and quadratic forms where $w$ is the local expression of $\omega$. But
the linear form $w$ is nowhere vanishing because the cokernel
sheaves are supported in $M_1$. So $\omega$ induces a morphism
\[
M_1\xrightarrow{\omega}\P V^\ast = S,
\]
which is locally given by $m\mapsto\langle w(m) \rangle$,
where $w(m)$ is the equation of the plane which contains the support of
$(\bku_1)_m$. It is now clear from this description that $\bku_1$ is also a
relative sheaf. Therefore there exists a unique
morphism
\[
\xymatrix{M_1 \ar[rr]^f\ar[dr]& & M_S \ar[dl]\\
 & S &  }
\]
with $\bku_1\cong f_{\,\P V\,}^\ast\bkv_S$. Conversely, $\bkv_S$ is a flat
family of stable sheaves with Hilbert polynomial $3m+1$ on $M_S\times \P V$
and thus defines a morphism $M_S\xrightarrow{g} M$ with ${\bkv}_S\cong
g_{\,\P V\,}^\ast\bku$. Now $g$ factors through $M_1$ and hence
$\bkv_S\cong
g_{\,\P V\,}^\ast\bku_1.$ It follows that $f$ and $g$ are inverse to each other.
\end{proof}
\vskip10mm

\section{Deformations and tangent spaces}

In this section we show that the elements of the intersection of the
two components $M_0$ and $M_1$ are
exactly the singular planar sheaves in $M_1$, i.e.~the planar sheaves which
are not locally free on their support. Furthermore, the two components, defined
as reduced subschemes, meet transversally.

Let $M'_1\subset M_1$ denote the locus of singular
sheaves. By lemma \ref{singpts} and proposition \ref{resprop} any
$\kf$ in $M'_1$ has a resolution
\[
0\to 2\ko_{\P_3}(-3)\xrightarrow{B_0} \ko_{\P_3}(-1)\oplus
3\ko_{\P_3}(-2)\xrightarrow{A_0}\ko_{\P_3}\oplus\ko_{\P_3}(-1)\to\kf\to 0
\]
with
\begin{eqnarray*}
  B_0=\begin{pmatrix} -q_1 & -l_1 & w & 0\\
                      -q_2 & -l_2 & 0 & w
      \end{pmatrix}
  & \text{ and } & A_0=\begin{pmatrix}
          w & 0\\
          0 & w\\
          q_1 & l_1\\
          q_2 & l_2
          \end{pmatrix}
\end{eqnarray*}
and $q_i=a_i\,l_1+b_i\,l_2$ for suitable linear forms $a_i,\, b_i$, $i=1,2$.

Using the same method as in the proof of \ref{locpres} one can prove that for 
any flat deformation
$\bkf\in\coh(V\times\P_3)$ of $\kf=\bkf_0$ over an open neighbourhood
$V\subset\A^1$ of $0$ there is an open neighbourhood $U\subset V$ over which  
$\bkf$ can be represented by a resolution
\[
0\to 2 \tko_U(-3)\xrightarrow{B_0+tB_1}\tko_U(-1)\oplus
3 \tko_U(-2)\xrightarrow{A_0+tA_1}\tko_U\oplus\tko_U(-1)\to\bkf\to 0,
\]
where $\tko_U(d):=\ko_U\boxtimes \ko_{\P_3}(d)$ and where $B_1$ and $A_1$ are
matrices which may depend on the parameter $t$ of $\A^1$. (To see this, note 
that $H^1(\P_3,\kf)=H^1(\P_3,\kf(1))=0$ and that then also
the direct images $R^1p_\ast\bkf, R^1p_\ast\bkf(1)$ vanish on a neighbourhood $U$ of 0 by 
the base change theorem. Assuming $U$ to be affine after shrinking, also
 $H^1(U\times \P_3,\bkf)=H^1(U\times \P_3,\bkf(1))=0.$ This implies that the 
homomorphism $A_0$ can be lifted and that its lifting is still surjective aftereventually shrinking $U$. Then use the same argument for the kernels.)

\vskip5mm

\begin{lemma}\label{deflem} $M'_1\subset M_0\cap M_1$.
\end{lemma}
\begin{proof} Let $[\kf]\in M'_1$ be arbitrary with corresponding
  matrices $A_0$ and $B_0$ as above. The ideal of the
  support $Z(\kf)$ is given by $(q_1l_2-q_2l_1,w)$.

  We construct explicitly a deformation $\bkf$ of the sheaf $\kf$ over
  $\A^1$ with $\bkf_0\cong \kf$ and $\bkf_t\in M_0\smallsetminus M_1$
  for $t\neq 0$ by choosing the following first order deformation matrices:
  \begin{eqnarray*}
  B_1=\begin{pmatrix} 0 & 0 & a_1 & b_1\\
                      0 & 0 & a_2 & b_2
      \end{pmatrix}
  & \text{ and } & A_1=\begin{pmatrix}
          a_1+b_2 & 1\\
          b_1a_2-b_2a_1 & 0\\
          0 & 0\\
          0 & 0
          \end{pmatrix}
\end{eqnarray*}
Consider the diagram:
\[\xymatrix{
      0\ar[r]\ar@{}[d]_-{t\neq0} & 2\,\tko(-3) \ar[rr]^-{
        B_t=B_0+t\,B_1}\ar@{}[rrd]|{\circlearrowright}\ar[d]_-{\approx}
        ^-{T_2} & &\tko(-1)\oplus 3\,\tko(-2) \ar[rr]^-{
          A_t=A_0+t\,A_1}\ar@{}[rrd]|{\circlearrowright}\ar[d]_-{T_1}^-{\approx} & &
    \tko\oplus
    \tko(-1) \ar[d]_-{T_0}^-{\approx}\ar@{}[rd]|{\circlearrowright}
    \ar[r] & \bkf  \ar@{..>}[d]^-{\approx}\ar[r] & 0 \\
      0\ar@{..>}[r] & 2\,\tko(-3) \ar@{..>}[rr]_-{\widetilde{B}_t} & &
      \tko(-1)\oplus
      3\,\tko(-2) \ar@{..>}[rr]_-{\widetilde{A}_t} &  & \tko\oplus
    \tko(-1)\ar@{..>}[r] & \bkg \ar@{..>}[r] & 0 }\]
  It is straightforward to check that the first row in the diagram is
  exact {\sl for all} $t\in\A_1$. Thus, the cokernel of $A_t$ is
  our $\A_1$-flat deformation $\bkf$.
  For $t\neq 0$, we choose the transformation matrices $T_i$ in the diagram as:
  \[
    T_{0,t}= \begin{pmatrix}
             1 & 0\\
             -a_1-b_2-t^{-1}w & -t^{-2}
            \end{pmatrix},\quad
    T_{1,t}= \begin{pmatrix}
             -t & 0 & 0 & 0\\
             -w & t & 0 & 0\\
             -l_1 & 0 & t & 0\\
             -l_2 & 0 & 0 & t
            \end{pmatrix}\quad\text{ and }\quad
    T_{2,t}= \begin{pmatrix} t^{-1} & 0\\
              0 & t^{-1}
              \end{pmatrix}
  \]
  Using those, we obtain the matrices in the second row:
\begin{eqnarray*}
  \widetilde{B}_t=\begin{pmatrix} 0 & -l_1 & ta_1+w & tb_1\\
                      0 & -l_2 & ta_2 & tb_2+w
      \end{pmatrix}
  & , & {\widetilde A}_t=\begin{pmatrix}
          0 & 1\\
          t^2b_1a_2-t^2b_2a_1-ta_1w-tb_2w-w^2 & 0\\
          tq_1-ta_1l_1-tb_2l_1-l_1w & 0\\
          tq_2-ta_1l_2-tb_2l_2-l_2w & 0
          \end{pmatrix}
  \end{eqnarray*}
  Again, one can verify the exactness of that row for $t\neq 0$.
  Consequently, we get as cokernel of ${\widetilde A}_t$ a family $\bkg\in\coh(\P_3\times\A_1)$ with
  $\bkf_t\cong\bkg_t=\ko_{C_t},\;t\neq 0$ where the curves $C_t$ are
  given by the ideal
  $(\,t^2b_1a_2-t^2b_2a_1-ta_1w-tb_2w-w^2,\;tq_1-ta_1l_1-tb_2l_1-l_1w,\;tq_2-ta_1l_2-tb_2l_2-l_2w\,)$.
  Therefore the isomorphism classes $[\bkf_t]$ are elements of $M\smallsetminus M_1\cong H_0\smallsetminus H_1$
  for non-vanishing $t$.\end{proof}
\vskip5mm

In order to show that in fact $M'_1=M_0\cap M_1$, we determine the dimension of
the tangent spaces at different points of $M$. This will also be used for the
proof of the transversality.
\vskip5mm

\begin{proposition}\label{tangdim}  Let $\kf$ be any sheaf in $M$.
  \begin{enumerate}
   \item [(i)] If $\kf\in M\smallsetminus M_1$, then $\dim T_{\kf}M=12$.
   \item [(ii)] If $\kf\in M_1\smallsetminus M'_1$, then $\dim T_{\kf}M=13$.
   \item [(iii)] If $\kf\in M'_1$, then $\dim T_{\kf}M=14$.
  \end{enumerate}
\end{proposition}

{\sl Preparations.} The parameter space $X$ defined in section 5 is reduced.
Therefore, also $M$ is reduced as a geometric quotient. $X$ consists of matrix
pairs $(B,A)$ with $B\circ A=0$ and $\text{rk}(\widetilde{B})=2$.

Let $p:=(B_0,A_0)\in X$ and denote by $F_p=G . p\,$ the $G$--orbit through $p$
which coincides with the fiber of $\kf:=\mu(p)$ under the map $\mu : X \to M$.
Standard arguments show that
\begin{eqnarray*}
T_pX & = & \{\,(B_1,A_1)\in X\; |\; B_1A_0+B_0A_1=0\,\}, \\
T_pF_p & = & \{\,(R B_0-B_0 S, SA_0-A_0T)\;|\;
(R,S,T)\in\text{Lie}(G)\,\}.
\end{eqnarray*}
Every $[\kg]\in\ext^1(\kf,\kf)$ has a resolution
$$\begin{CD}
0 @>>> \kl_2\oplus\kl_2 @>{{\tiny\text{$\begin{pmatrix} B_0 & 0 \\ B_1 &
    B_0\end{pmatrix}$}}}>> \kl_1\oplus\kl_1 @>{{\tiny\text{$\begin{pmatrix} A_0 & 0 \\ A_1 &
    A_0\end{pmatrix}$}}}>> \kl_0\oplus\kl_0 @>>> \kg @>>> 0,
\end{CD}$$
obtained by adding two copies of the resolution \ref{resprop} of $\kf$.
Moreover,
\begin{itemize}\itemsep4pt
\item ${\tiny\text{$\begin{pmatrix} B_0 & 0 \\ B_1 &
    B_0\end{pmatrix}$ $\begin{pmatrix} A_0 & 0 \\ A_1 &
    A_0\end{pmatrix}$}}=0$ is equivalent to $B_1 A_0 + A_0 B_1 =0$.
\item $\kg\cong\kf\oplus\kf$ if and only if $A_1=S A_0-A_0 T$ and $B_1=R
B_0-B_0 T$ for some matrices $S,T,R$.
\end{itemize}
Therefore we get an exact sequence
$0\to T_pF_p\longrightarrow T_pX \longrightarrow \ext^1(\kf,\kf) \to 0$,
where the second map is given by $(B_1,A_1)\mapsto [\kg]$. Since the
Zariski tangent space $T_\kf M$ is isomorphic to the extension group
$\ext^1(\kf,\kf)$, cf.~proposition \ref{tangen}, we obtain
$$T_{\kf}M\cong T_pX/T_pF_p.\eqno (\ast)$$

\begin{proof} Part (i) is an immediate consequence of lemma \ref{spcub}:
$M\smallsetminus M_1\cong H_0\smallsetminus H_1$ is smooth of
dimension $12$. For (ii) and (iii) consider $X_1=\mu^{-1}(M_1)=\{\,(B,A)\in
X\;|\;\lambda(A)=0\,\}$,
where $\lambda(A)$ is the degree $0$ entry of $A$. Let $X_1'$ denote
the preimage of the subset $\mu^{-1}(M_1')$ of singular sheaves in $M_1$.
Clearly, $T_pX_1\subset T_pX$ for $p\in X_1$.

(1) For $p=(B_0,A_0)\in X_1\smallsetminus X_1'$, we have $T_pX_1=T_pX$.
To see this, let $(B_1,A_1)\in T_pX$. The equation $B_1 A_0 + B_0 A_1 =0$ can
be written more explicitly as
$$\begin{pmatrix} r_1 & x_1 & x_3 & x_5\\r_2 & x_2 & x_4 &
  x_6\end{pmatrix}
\begin{pmatrix} w & 0 \\ 0 & w\\q_1 & l_1 \\ q_2 & l_2\end{pmatrix}
+
\begin{pmatrix} -q_1 & -l_1 & w & 0\\-q_2 & -l_2 & 0 & w \end{pmatrix}
\begin{pmatrix}   u  & \lambda\\ h & v\\ h_1 & v_1\\ h_2 & v_2\end{pmatrix}
=0.$$
The two entries $(1,2)$ and $(2,2)$ of this matrix modulo the three linear forms
$w,l_1$ and $l_2$ yield the equations
$$\overline{q}_1 \lambda=0\quad\text{and}\quad\overline{q}_2 \lambda=0.$$
Since $p\not\in X_1'$, the quadratic forms $q_1$ and $q_2$ are not
contained in the ideal $(l_1,l_2,w)$ and therefore
$\overline{q}_1$ and $\overline{q}_2$ cannot both be equal to zero.
Hence $\lambda=0$ and $(B_1,A_1)\in T_pX_1$.

(2) On the other hand, if $p\in X_1'$ there exists in addition the tangent vector
$$(B_1,A_1)=\left(\begin{pmatrix} 0 & 0 & a_1 & b_1\\
                      0 & 0 & a_2 & b_2
      \end{pmatrix}\,,\,\begin{pmatrix}
          a_1+b_2 & 1\\
          b_1a_2-b_2a_1 & 0\\
          0 & 0\\
          0 & 0
          \end{pmatrix}\,\right)$$
from lemma \ref{deflem} which is {\sl not} contained in $T_pX_1$.
A direct computation shows that $T_pX_1$ and the vector $(B_1,A_1)$ span
$T_pX$. Therefore $\dim T_pX=\dim T_pX_1+1$ in this case.
Since $F_p\subset X_1$ for all $p\in X_1$, we obtain from $(\ast)$ that
$$T_\kf M/T_\kf M_1\cong T_p X/T_p X_1.$$
According to proposition \ref{m1}, $M_1$ is smooth of dimension $13,$
and we get
$$\dim T_{\kf} M= \dim
T_\kf M_1+\dim T_p X/T_p X_1=\left\{ \begin{array}{l@{\,\text{ if }}l} 14 &
      \kf\in M_1' \\  13 & \kf\in M_1\smallsetminus
      M_1'.\end{array}\right.$$
\end{proof}

\begin{corollary}\label{intsec} The intersection $M_0\cap M_1$ coincides
  with the set $M'_1$ of singular planar sheaves.
\end{corollary}

\begin{proof} By lemma \ref{deflem}, $M_1'\subset M_0\cap M_1$. Then
  $\dim T_\kf M > \dim M =13$ for any $\kf$ in
  $M_0\cap M_1\smallsetminus M_1'$. By proposition \ref{tangdim} (ii), this is not possible.
\end{proof}

\vskip5mm

\begin{proposition}\label{transv} For any point $p\in M_0\cap M_1$,
  \begin{enumerate}
  \item [(a)]   $T_pM_0+T_pM_1=T_pM$
  \item [(b)]  $T_p(M_0 \cap M_1)=T_pM_0\cap T_pM_1$
  \end{enumerate}
\end{proposition}

\begin{proof}We have $T_pM_0+T_pM_1\subset T_pM$ and therefore $13\leq
  \dim (T_pM_0+T_pM_1)\leq\dim T_pM=14$ by proposition \ref{tangdim}. If $\dim
  (T_pM_0+T_pM_1)=13$, then $T_pM_0\subset T_pM_1$ since $\dim
  T_p M_1=13$ because of the smoothness of $M_1$.
  But this is not possible because there are deformation paths
  from $p$ to points $q\in M_0\smallsetminus M_1$ as we saw in the proof of
  lemma \ref{deflem}. Part (b) is a general fact for intersections of subvarieties.
\end{proof}

{\bf Remark}: We will show in section 8 that
$M_0\cap M_1$ is smooth of dimension 11. Then it follows already from (a) and
(b) that $M_0$ is also smooth along $M_0\cap M_1$ and that the components
$M_0$ and $M_1$ intersect transversally.

\vskip10mm

\section{The component $M_0$ and transversality}

We will now prove that the isomorphism $M_0\smallsetminus M_1\cong
H_0\smallsetminus H_1$ from lemma \ref{spcub} extends to an
isomorphism $M_0\cong H_0$ and
that moreover $M_0$ is isomorphic to a blowup of the smooth variety
$N$ of nets of
quadrics as described in \cite{E/S}, see \ref{nets} below. Each of these
two statements implies that $M_0$ is smooth and rational of
dimension $12$ due to the
results in \cite{P/S} and \cite{E/S}. The blowup structure of $H_0$ has been
mentioned without proof in \cite{E/S}. We include a nice explicit description
for $M_0$ using the deformation matrices of section 7.

\begin{proposition}\label{fitting} The component $M_0$ is isomorphic
  to the component $H_0$ of the Hilbert scheme $\hilb_{3m+1}(\P_3)$.
\end{proposition}

\begin{proof}
We extend the isomorphism $f' : M_0\smallsetminus M_1\to
 H_0\smallsetminus H_1$ using Fitting ideals. The Fitting ideal of a sheaf $\kf\in M$ is the ideal
 $\text{Fitt}(\kf)$ generated by the
 $2$--minors of the presentation matrix $A$ in the resolution of
 $\kf$. Recall from proposition \ref{resprop} that $A$ has
 the form
\[
\left(
\begin{array}{cc}
w & 0\\
0 & w\\
q_1 & l_1\\
q_2 & l_2
\end{array}
\right)\quad\text{ or }\quad
\left(
  \begin{array}{cc}
0 & 1\\
q_1 & 0\\
q_2 & 0\\
q_3 & 0
  \end{array}
\right)
\]
depending on whether the support $Z(\kf)$ is contained in a plane or not.
By corollary \ref{intsec}, the elements of $M_0\cap M_1$
are exactly the singular planar sheaves which are characterized by
matrices $A$ of the first type with the additional condition
$(q_1, q_2)\subset (l_1, l_2)$. Therefore we obtain
\begin{eqnarray*}
  \text{Fitt}(\kf) & = & (q_1,q_2,q_3)\quad\text{if }\kf\in M_0
  \smallsetminus M_1\\
  \text{Fitt}(\kf) & = & (w^2, wl_1, wl_2, l_1q_2-l_2q_1)\quad\text{if
    }\kf\in M_0 \cap M_1
\end{eqnarray*}
Note, that on the other hand $(w^2, wl_1, wl_2, l_1q_2-l_2q_1)$ is
exactly the normal form for ideals of cubics in $H_0\cap H_1$, see
\cite{P/S}. Therefore we get a surjective, set-theoretical map
$f : M_0\longrightarrow H_0,\;
\kf\mapsto\ko_{Z(\text{Fitt}(\kf))}$
extending the isomorphism $f'$. It is indeed bijective, the
injectivity on $M_0\cap M_1$ is easy to check. Applying
$\ko_{Z(\text{Fitt}(\,\cdot\,))}$ to each fiber $\kf=\bku_{\kf}$ of the
universal sheaf $\bku|_{M_0\times\P_3}$, we obtain a family
$\bkg\in\coh(H_0\times \P_3)$.
It is straightforward to see that $\bkg$ is $H_0$--flat with constant Hilbert
polynomial $3m+1$ along the fibers. The morphism induced by $\bkg\in
\km_{3m+1}(\P_3)(H_0)$ is exactly the set-theoretical map
$f$.
Since $H_0$ is smooth the bijective morphism $f$ is an
isomorphism by Zariski's main theorem.
\end{proof}
\vskip5mm
{\bf Remark:} The above Fitting map cannot be extended to a morphism
$M_{3m+1}(\P_3)\to \hilb_{3m+1}(\P_3)$, because the family of Fitting ideals of
all the matrices $A$ is not flat over $M$.
\vskip5mm

\begin{sub}\label{nets}\rm{\bf Nets of quadrics.}\\
Let $V\cong k^4$ and let $N$ denote the geometric quotient
\[
\hom(k^2, k^3\otimes V^\ast)^s/\GL_2(k)\times GL_3(k)
\]
as described in \cite{E/S}. An adhoc definition of the stability of a matrix
\[
Q=\left(
  \begin{array}{ccc}
x_1 & x_2 & x_3\\
y_1 & y_2 & y_3
  \end{array}\right)
\]
in $\hom(k^2, k^3\otimes V^\ast)$ is that $Q$ is not equivalent to a matrix
with two zeros in a row or in a column. The scheme $N$ is projective, smooth,
rational of dimension 12. Its tangent space at $[Q]$ is isomorphic to
$$
T_{[Q]}N\cong\hom(k^2, k^3\otimes V^\ast)/F_Q \qquad \text{with}\quad
F_Q=\{SQ-QR\, |\quad S\in\mathfrak{gl}_2, R\in\mathfrak{gl}_3\},
$$
as in the proof of proposition \ref{tangdim}. Let $N_1\subset N$ be the
subvariety whose points are of the form
\[
\left[
  \begin{array}{ccc}
l_1 & w & 0\\
l_2 & 0 & w
  \end{array}\right]
\]
with linear independent forms $w, l_1, l_2 \in V^\ast$. $N_1$ is the
subvariety of classes of matrices, whose $3$ minors have a common linear
factor and thus do not determine a cubic space curve. On the other hand
\[
N\smallsetminus N_1\cong H_0\smallsetminus H_1\cong M\smallsetminus M_1
\]
because the 3 quadrics of a matrix in $N\smallsetminus N_1$ generate the ideal
of a cubic space curve in $\P_3=\P V$, see \cite{E/S}.

Let $q\in N_1$ be defined by the matrix $Q=\left(\begin{smallmatrix}l_1 & w &
  0\\l_2
  & 0 & w\end{smallmatrix}\right)$ and let $l_3$ be a fourth linear form independent
  of $w, l_1, l_2$, and let
\[
T'_Q\subset \hom (k^2, k^3\otimes V^\ast)
\]
be the $5$--dimensional subspace spanned by the matrices
\[
\left(
  \begin{array}{c|c|c}
\rho l_3 & \lambda_1 l_1+\lambda_2l_2+\lambda_3l_3 & 0\\\hline
\sigma l_3 & 0 & \lambda_1l_1+\lambda_2 l_2+\lambda_3l_3
  \end{array}\right)
\]
with $\rho, \sigma, \lambda_1, \lambda_2, \lambda_3\in k$. We have the
\end{sub}
\vskip5mm

\begin{lemma}\label{plnets}
  \begin{enumerate}
  \item [(i)]   $T'_Q \cap F_Q=0$
  \item [(ii)]  $T_qN_1= T'_Q\oplus F_Q/F_Q$
  \item [(iii)] $T_qN/T_qN_1\cong \hom(k^2, k^3\otimes V^\ast)/T'_Q\oplus F_Q$
  \end{enumerate}
\end{lemma}
\begin{proof} (i) follows by a direct computation in linear algebra. Because
  the matrices of $T'_q$ define a slice over $N_1$ by $Q+T'_Q$ and cover a
  neighbourhood of $q$ in $N_1$, we obtain that $T'_Q\oplus F_Q/F_Q\subset T_Q
  N_1\subset T'_Q\oplus F_Q/F_Q$. (iii) is a consequence of (ii).
\end{proof}

{\bf Remark:} The above lemma implies also that $N_1$ is smooth of dimension
$5$ because $\dim T'_Q=5$, $\dim F_Q=12$. It is actually isomorphic to the
tautological plane bundle over $\P_3^\ast$, in other words, to the
flag manifold $\F lag(1,3,V^\ast)$.
\vskip5mm

\begin{sub}\label{blm}\rm {\bf The morphism $M\xrightarrow{\rho} N$.}\\
  Recalling the definition of the parameter space $X=X_0 \cup X_1$ of $M$, as
  a space of pairs $(B,A)$ with
\[
B=\left(
  \begin{array}{c|ccc}
f & x_1 & x_2 & x_3\\
g & y_1 & y_2 & y_3
  \end{array}\right),\quad A=\left(
  \begin{array}{c|c}
z & \lambda\\\hline
q_1 & z_1 \\
q_2 & z_2\\
q_3 & z_3
  \end{array}\right),
\]
we know that the component $X_1$ is defined by the condition $\lambda=0$. It
follows that $M_1\subset M$ and then also $M_0\cap M_1\subset M_0$ is a
Cartier divisor in its reduced structure. Moreover, the submatrix
$\left(\begin{smallmatrix}x_1 & x_2 & x_3\\y_1
  & y_2 & y_3\end{smallmatrix}\right)$
of $B$ is stable, see last part of the proof of lemma \ref{lempara}. Thus
\[
(B, A)\mapsto \left(
  \begin{array}{ccc}
x_1 & x_2 & x_3\\
y_1 & y_2 & y_3
  \end{array}\right)
\]
is an equivariant morphism
$X\to\hom (k^2, k^3\otimes V^\ast)^s,$
which induces a morphism $M\xrightarrow{\rho} N$
of the geometric quotients. It maps $M_1$ and $M_0 \cap M_1$ onto $N_1$
as follows from the normal form of pairs $(B,A)\in X_1$, see
proposition \ref{resprop}. Moreover, the restriction
\[
M\smallsetminus M_1\xrightarrow[\approx]{\rho} N\smallsetminus N_1
\]
is an isomorphism, using the normal form with $\lambda=1$.
\end{sub}
\vskip5mm

\begin{proposition}\label{blup} $M_0\xrightarrow{\rho}N$ is
  isomorphic to the blowup $\Bl_{N_1} N$ of $N$ along $N_1$, and $M_0 \cap M_1$ is
  isomorphic to the projective normal bundle $\P(\kn_{N_1/N})$.
\end{proposition}

\begin{proof} a) Because $M_0 \cap M_1$ is a Cartier divisor in $M_0$ which is
  mapped onto $N_1$, there is the unique natural morphism
\[
\xymatrix{M_0\ar[rr]^\varphi\ar[dr]_\rho & & \Bl_{N_1}(N)\ar[dl]^\sigma\\
& N
}
\]
mapping $M_0 \cap M_1$ to the exceptional divisor $E\cong \P(\kn_{N_1/N})$, where we
use that $N_1$ is smooth. In our case $\varphi|_{M_0\smallsetminus
(M_0 \cap M_1)}$ is
already an isomorphism. We are now going to prove that the induced morphism
\[
M_0 \cap M_1\xrightarrow{\tau}\P(\kn_{N_1/N})
\]
is an isomorphism. For that it is sufficient, that the induced map
\[
(M_0 \cap M_1)_q\xrightarrow{\tau_q}\P(T_qN/T_qN_1)
\]
is bijective for the fibers over a point $q\in N_1$. Then $\varphi$ is
bijective and an isomorphism because $\Bl_{N_1}(N)$ is smooth.

b) The map $\tau_q$ is defined as follows. For any point $[B_0, A_0]\in
M_0 \cap M_1$ over $q$ choose a tangent vector
\[
\xi\in T_{[B_0, A_0]}M_0\smallsetminus T_{[B_0, A_0]} M_1.
\]
Its image in $T_qN\smallsetminus T_qN_1$ and then its image in $\P(T_q N/T_q
N_1)$ is the image $\tau_q([B_0, A_0])$. More explicitly, if
\[
q=\left[
  \begin{array}{ccc}
l_1 & w & 0\\
l_2 & 0 & w
  \end{array}\right]\;\in N_1,
\]
the points $[B_0, A_0]$ over $q$ are of the form
\[
B_0 = \left(
  \begin{array}{cccc}
-q_1 & -l_1 & w & 0\\
-q_2 & -l_2 & 0 & w
  \end{array}\right),\quad A_0=\left(
  \begin{array}{cc}
w & 0\\
0 & w\\
q_1 & l_1\\
q_2 & l_2
  \end{array}\right)
\]
with $q_1, q_2\in(l_1, l_2)$, see section 7. Then
\[
q_1=a_1l_1+b_1l_2\quad \text{ and }\quad q_2=a_2l_1+b_2l_2
\]
as in section 7, and for the tangent vector $\xi$ we can choose the image of
the pair
\[
B_1=\left(
  \begin{array}{cccc}
0 & 0 & a_1 & b_1\\
0 & 0 & a_2 & b_2
  \end{array}\right)\ ,\quad A_1=\left(
  \begin{array}{cc}
a_1+b_2 & 1\\
b_1a_2-b_2a_1 & 0\\
0 & 0\\
0 & 0
  \end{array}\right)
\]
as in the proof of lemma \ref{deflem}. The image of this tangent vector in
$\P(T_q N/T_q N_1)$ is simply the class of the matrix
$\left(
\begin{smallmatrix} 0 & a_1 & b_1\\0 & a_2 & b_2
\end{smallmatrix}\right).$
Hence, the map $\tau_q$ is given in notation of classes by
\[
\left[
\left(
  \begin{array}{cccc}
-q_1 & -l_1 & w & 0\\
-q_2 & -l_2 & 0 & w
  \end{array}\right),\ \left(
  \begin{array}{cc}
w & 0\\
0 & w\\
q_1 & l_1\\
q_2 & l_2
  \end{array}\right)\right]\mapsto \langle\overline{\left(
  \begin{array}{ccc}
0 & a_1 & b_1\\
0 & a_2 & b_2
  \end{array}\right)}\rangle\ .
\]
It is now an easy exercise to verify that this map is well-defined and
bijective. This completes the proof of proposition \ref{blup}.
\end{proof}
\vskip5mm

\begin{corollary} $M_0$ is rational and smooth of dimension 12 and
  $M_0\cap M_1$ is a smooth rational divisor in $M_0$.
\end{corollary}
\vskip5mm

\begin{corollary} The components $M_0$ and $M_1$ intersect
  transversally in $M_0\cap M_1$ and for any $p\in M_0\cap M_1$
\[
T_pM_0+T_pM_1=T_pM
\]
has dimension $14$, whereas $\dim T_pM_0\cap T_pM_1=11$.
\end{corollary}
\vskip5mm

\section{The Chow groups of the moduli space}

Let $k=\C$ and $V\cong k^4$.
Using the descriptions of $M_1$ as the relative universal cubic
$M_1=\kz\rightarrow\P V^\ast$ and of $M_0$ as the blow up $\text{Bl}_{N_1}(N)$
with exceptional divisor $M_0\cap
M_1\cong\P(\kn_{N_1/N})$, we determine now the Betti numbers of the
Chow groups of the moduli space $M=M_{3m+1}(\P_3)$. Key ingredient is
the following result on the Chow groups $A_\ast(N)$ of the space $N$ of nets of
quadrics, see \cite{E/S}:

\begin{proposition}[Ellingsrud, Str{\o}mme]
\label{Ell/Stro}
The Betti numbers are

$$\begin{array}{|c||c|c|c|c|c|c|c|c|c|c|c|c|c|}
   \hline
   i & 0 & 1 & 2 & 3 & 4 & 5 & 6 & 7 & 8 & 9 & 10 & 11 & 12\\
   \hline\hline
   b_i(N) & 1 & 1 & 3 & 4 & 7 & 8 & 10 & 8 & 7 & 4 & 3 & 1 & 1\\
   \hline
  \end{array}$$

In particular, the topological Euler characteristic is $e(N)=58$.
\end{proposition}

\begin{theorem}\label{chowgr}
The Chow groups of $M_{3m+1}(\P_3)$, of its components $M_0, M_1$ and
of the intersection $M_0\cap M_1$ are free and have the following topological
Betti numbers $b_i$:

$$\begin{array}{|c||c|c|c|c|c|c|c|c|c|c|c|c|c|c||c|}
   \hline
   i & 0 & 1 & 2 & 3 & 4 & 5 & 6 & 7 & 8 & 9 & 10 & 11 & 12 & 13 &
    e(\,\cdot\,)\\
   \hline\hline
 b_i(M_0\cap M_1)  & 1 & 3 & 6 & 9 & 11 & 12 & 12 &
 11 & 9 & 6 & 3 & 1 & - & - & 84\\
 \hline
 b_i(M_0) & 1 & 2 & 6 & 10 & 16 & 19 & 22 & 19 & 16
 & 10 & 6 & 2 & 1 & - & 130\\
 \hline
   b_i(M_1)  & 1 & 3 & 6 & 9 & 11 & 12 & 12 & 12 & 12
   & 11 & 9 & 6 & 3 & 1 & 108 \\
 \hline\hline
  b_i(M_{3m+1}(\P_3)) & 1 & 2 & 6 & 10 & 16 & 19 & 22 & 20 & 19
   & 15 & 12 & 7 & 4 & 1 & 154\\
 \hline
  \end{array}\,.$$

Furthermore, the Chow ring $A^\ast(M_1)$ is isomorphic to $\Z[s,t,u]/\fraca$,
where the ideal $\fraca$ is given by
$$\fraca=(s^4,t^3-st^2+s^2t-s^3,u^9+(10\,s-3\,t)u^8+(55\,s^2-30\,st+9\,t^2)u^7+(220\,s^3-165\,s^2t+90\,st^2)u^6+$$
$$+(495\,s^2t^2-660\,s^3t)u^5+1980\,s^3t^2u^4).$$
\end{theorem}

\begin{proof}
Since $N_1\cong \F lag(1,3,V^\ast)$ we obtain immediately
$$A^\ast(N_1)\cong A^\ast(\P V)[t]/(t^3+c_1(\kq)t^2+c_2(\kq)t+c_3(\kq))\cong\Z[s,t]/(s^4,s^3+st^2+s^2t+s^3),$$
where $\kq$ is the tautological quotient in the sequence
$0\to \ko_{\P V}(-1)\to\ko_{\P V}\otimes V^\ast\to\kq\to 0.$
We get the Betti numbers for $M_0\cap M_1$ using
$$A_i(M_0\cap M_1)\cong A_i(\P(\kn_{N_1/N}))\cong\bigoplus_{k=0}^6
A_{i-6+k}(N_1)$$ and the presentation of the ring $A^\ast(N_1)$. But
then it is straightforward to compute the numbers for $M_0$ from
\[\xymatrix{0\ar[r] & A_i(N_1)\ar[r]\ar[d]^{\approx} &
  A_i(\P(\kn_{N_1/N}))\oplus A_i(N) \ar[r]\ar[d]^{\approx} &
  A_i(\text{Bl}_{N_1}(N))\ar[r]\ar[d]^{\approx} & 0\\
0\ar[r] & A_i(N_1)\ar[r] & A_i(M_0\cap M_1)\oplus A_i(N)\ar[r] &
A_i(M_0)\ar[r] & 0}\]
since we know $b_i(N)$, $i=1,\ldots,12$ from proposition
\ref{Ell/Stro}.

Now recall the construction of the relative
universal cubic $\kz\rightarrow\P V^\ast$ parametrizing triples
$(p,C,H)$ of points $p$ on cubic curves $C$ contained in planes
$H\subset\P V^\ast$:
Consider the exact sequence
$$0\to\kh\longrightarrow \ko_{\P V^\ast}\otimes
V\xrightarrow{eval}\ko_{\P V^\ast}(1)
\to 0.$$
The plane bundle $\P(\kh)\stackrel{p}{\longrightarrow}\P V^\ast$ comes with a surjective map
$p^\ast\kh^\ast\to\ko_{\kh}(1)\to 0$ which in turn induces
$$0\to\kk\longrightarrow p^\ast S^3\kh^\ast
\longrightarrow\ko_{\kh}
(3)\to 0.$$
Then $\kz=\P(\kk)\rightarrow\P(\kh)\rightarrow\P V^\ast$ is the
desired incidence variety.

Using this description of $\kz$, we compute the Chow ring of
$M_1\cong\kz$ as follows: The total Chern class of $\kh$ is $c(\kh)=1-s+s^2-s^3$ and
consequently
$$A^\ast(\P(\kh))=\Z[s,t]/(s^4,t^3-st^2+s^2t-s^3).$$ One
can easily check that
$c(S^3\kh^\ast)=1+10\,s+55\,s^{\,2}+220\,s^{\,3}$.
Then $c(\kk)=(1+10\,s+55\,s^2+220\,s^3)(1+3t)^{-1}$ implies
$$A^\ast(M_1)\cong A^\ast(\P(\kh))[u]/(u^9+c_1(\kk)u^8+\cdots
+c_9(\kk))\cong \Z[s,t,u]/(s^4,t^3-st^2+s^2t-s^3,f)\,,$$
where
$f=u^9+(10\,s-3\,t)u^8+(55\,s^2-30\,st+9\,t^2)u^7+(220\,s^3-165\,s^2t+90\,st^2)u^6+
(495\,s^2t^2-660\,s^3t)u^5+1980\,s^3t^2u^4$.
This allows us also to determine the Betti numbers of $M_1$.
Finally, since $k\cong\C$, we get the following diagram
\[\xymatrix{
   H_{2i+1}(M) \ar[r] & H_{2i}(M_0\cap M_1) \ar[r] & H_{2i}(M_0)\oplus
H_{2i}(M_1) \ar[r] & H_{2i}(M) \ar[r] & H_{2i-1}(M_0\cap M_1)\\
0\ar[u]_{\approx}\ar[r] & A_i(M_0\cap M_1)
\ar[r]\ar[u]_{\approx}^-{cl} & A_i(M_0) \oplus
A_i(M_1)\ar[r]\ar[u]_-{\approx}^-{cl} &
A_i(M)\ar[r]\ar[u]_-{\approx}^-{cl} & 0\, ,\ar[u]_{\approx}}
\]
where we use homology with locally compact supports. Note
that in our case the cycle maps $cl$ are isomorphisms.
Using this Mayer--Vietoris sequence, we obtain immediately the Betti numbers of the moduli space $M=M_{3m+1}(\P_3)$.
\end{proof}

\end{document}